\documentclass[16pt]{amsart}
\usepackage{amsmath,amstext,amsthm}
\usepackage{amssymb}
\usepackage{amsfonts}
\address{G\'{e}rard Besson \\Institut Fourier\\ UMR 5582 CNRS \\
BP 74 38402 St Martin d'H\`{e}res \\ France }\usepackage[dvips]{graphics}
\newtheorem{theorem}{Theorem}[section]
\newtheorem{rem}[theorem]{Remark}
\newtheorem{lemma}[theorem]{Lemma}
\newtheorem{prop}[theorem]{Proposition}
\newtheorem{definition}[theorem]{Definition}
\newtheorem{cor}[theorem]{Corollary}

\numberwithin{equation}{section}

\begin{document}
\title{Rigidity of amalgamated product in negative curvature}
\author{G\'{e}rard Besson, Gilles Courtois, Sylvain Gallot}
\email{G.Besson@ujf-grenoble.fr}
\address{Gilles Courtois \\ CMLS \\ \'{E}cole Polytechnique\\
UMR 7640 CNRS\\ 91128 Palaiseau \\ France}
\email{courtois@math.polytechnique.fr}
\address{Sylvain Gallot \\Institut Fourier\\ UMR 5582 CNRS \\
BP 74 38402 St Martin d'H\`{e}res \\ France }
\email{sylvestre.gallot@ujf-grenoble.fr}
\date{\today}
\maketitle

\large

\begin{abstract}
Let $\Gamma$ be the fundamental group of a compact riemannian manifold $X$
of sectional curvature $K\leq -1$ and dimension $n\geq 3$.
We suppose that $\Gamma = A\ast_{C}B$ is the free product of its subgroups
$A$ and $B$ over the amalgamated subgroup $C$. We prove that the critical exponent
$\delta(C)$ of $C$ satisfies $\delta(C) \geq n-2$. The equality happens if and only if 
there exist an embedded compact
 hypersurface $Y \subset X$, totally geodesic, of constant sectional curvature $-1$,
whose fundamental group is $C$ and which separates $X$ in two connected components
whose fundamental groups are $A$ and $B$. 
Similar results hold if $\Gamma$ is an HNN extension, or more generally if $\Gamma$
acts on a simplicial tree without fixed point.
 
\end{abstract}

\section{Introduction}

In \cite{sh}, Y. Shalom proved the following theorem which says that for every lattice
$\Gamma$ in the hyperbolic space and for any decomposition of $\Gamma$ as an amalgamated
product $\Gamma = A \ast_{C} B$, the group $C$ has to be ``big''. In order to measure how ``big''
$C$ is, let us define the critical exponent of a discrete group $C$ acting on a Cartan Hadamard manifold
by

$$
\delta(C) = inf \{ s >0 \quad | \quad \Sigma_{\gamma \in \Gamma} e^{-sd(\gamma x,x)} < +\infty \}.
$$

\begin{theorem}
Let $\Gamma $ be a lattice in $PO(n,1)$. Assume 
that $\Gamma$ is an amalgamated product of 
its subgroups $A$ and $B$ over $C$.
Then, the critical exponent $\delta(C)$ of $C$ satisfy
$\delta(C) \geq n-2$.
\end{theorem}

An example is given by any $n$-dimensional hyperbolic manifold $X$ which contains 
a compact separating connected 
totally geodesic hypersurface $Y$. The Van Kampen theorem
then says that the fundamental group $\Gamma$ of $X$ is isomorphic to the free product
of the fundamental groups of the two halves of $X-Y$ amalgamated over the fundamental
group $C$ of the incompressible hypersurface $Y$. 
Such examples do exist in dimension $3$ thanks to the W.Thurston's hyperbolization theorem.
In any dimension, A. Lubotsky showed that any standard arithmetic lattice of $PO(n,1)$ has a finite cover whose
fundamental group is an amalgamated product, cf. \cite{lu}.
In fact, A. Lubotsky proved that any standard arithmetic lattice $\Gamma$ has a finite index subgroup $\Gamma_{0}$
which is mapped onto a nonabelian free group. A nonabelian free group can be written in infinitely many ways as
an amalgamated product, so one get infinitely many decomposition of $\Gamma_{0}$ as an amalgamated product by pulling
back the amalgamated decomposition of the nonabelian free group. 

In these cases there is equality in theorem 1.1, ie. $\delta(C) = n-2$ where $C$ is 
the fundamental group of $Y$, and Y. Shalom suggested in \cite{sh} 
that the equality case in the theorem 1.1 happens only in that case.  

The aim of this paper is to show that the theorem 1.1 still holds 
when $\Gamma$ is the fundamental group of a compact riemannian
manifold of variable sectional curvature less than or equal to $-1$,
and characterize the equality case.

\begin{theorem}
Let $X$ be an $n$-dimensional compact riemanniann manifold of 
sectional curvature $K\leq -1$.
We assume that the fundamental group $\Gamma$ of $X$ is 
an amalgamated product of 
its subgroups $A$ and $B$ over $C$ and that neither $A$ nor $B$ equals
 $\Gamma$.
Then, the critical exponent $\delta(C)$ of $C$ satisfy
$\delta(C) \geq n-2$.
Equality $\delta(C) = n-2$ happens if and only if 
$C$ cocompactly preserves a totally geodesic isometrically embedded
 copy $\mathbb H^{n-1}$
of the hyperbolic space of dimension $n-1$.
Moreover, in the equality case, the hypersurface $Y^{n-1} := \mathbb H^{n-1}/C$
is embedded in $X$ and separates $X$ in two connected components 
whose fundamental groups are respectively $A$ and $B$.
\end{theorem}

\begin{rem}
(i) By the assumption on $A$ or $B$ not being equal to $\Gamma$ we exclude
the trivial decomposition  $\Gamma = \Gamma \ast_{C} C$ where $A=\Gamma$
and $C = B$ can
be an arbitrary subgroup of $\Gamma$, for example any cyclic subgroup, 
in which case the conclusion of theorem 1.2 fails. 
Also note that because of this assumption on $A$ and $B$, we have
$A\neq C$ and $B \neq C$ . 

(ii) Let us recall that standard
arithmetic lattices in $PO(n,1)$ have finite index subgroup with infinitely many
non equivalent decompositions as amalgamated products, cf. \cite{lu}. In fact, among these decompositions,
all but finitely many of them are such that $\delta(C) > n-2$. Indeed, by theorem 1.2, if
$\delta(C)=n-2$ then $C$ is the fundamental group of an embedded totally geodesic 
hypersurface in $X$, but there are only finitely many totally geodesic hypersurfaces by \cite{ze}.  
\end{rem}

When a group is an amalgamated product, it acts on a simplicial tree without fixed point
and theorem 1.1 is a particular case of the

\begin{theorem}
(\cite{sh}, theorem 1.6).  Let $\Gamma \subset SO(n,1)$ , $n\geq 3$, be a lattice
Suppose $\Gamma$ acts on a simplicial tree $T$ without fixed vertex. Then there is an edge
of $T$ whose stabilizer $C$ satisfies $\delta(C) \geq n-2$.
\end{theorem}

In the case $\Gamma$ is cocompact, the conclusion of theorem 1.4 holds for 
the stabilizer of any edge which separates
the tree $T$ in two unbounded components, and the proof of this is exactly the same as 
the proof of theorem 1.2.
 In particular, when the action of $\Gamma$
on $T$ is minimal, (ie. there is no proper subtree of $T$ invariant by $\Gamma$),
the conclusion of theorem 1.4 holds for every edge of $T$,  
in the variable curvature setting,
and we are able to handle the equality case.

\begin{theorem}
Let $\Gamma$ be the fundamental group of an $n$-dimensional compact riemannian manifold $X$
of sectional curvature less than or equal to $-1$. 
Suppose $\Gamma$ acts minimally on a simplicial tree $T$ without fixed point.
Then, the stabilizer $C$ of every edge of $T$ satisfies $\delta(C) \geq n-2$.
The equality $\delta(C) = n-2$ happens 
if and only if there exist a compact totally geodesic hypersurface
$Y\subset X$ with fundamental group $\pi_{1}(Y) = C$. Moreover, 
in that case, $Y$ with its induced metric has constant sectional curvature $-1$.
\end{theorem} 

\bigskip
Another interesting case contained in theorem 1.5 is the case of $HNN$ extension. 
Let us recall the definition of an $HNN$ extension. Let $A$ and $C$ be groups
and $ f_{1} :C \to A$, $f_{2}: C \to A$ two injective morphisms of $C$ into $A$.
The $HNN$ extension $A \ast_{C}$ is the group generated by $A$ and an element $t$
with the relations $t f_{1}(\gamma) t^{-1} = f_{2}(\gamma)$.
For example, let $X$ be a compact manifold containing a non separating compact 
incompressible hypersurface $Y\subset X$. Let $A$ be the fundamental group of the manifold 
with boundary $X - Y$ obtained by cuting $X$ along $Y$ and let $C$ be the fundamental
group of $Y$. The boundary of $X-Y$ consists in two connected components
$Y_{1}\subset X-Y$ and $Y_{2}\subset X-Y$ homeomorphic to $Y$.
By the incompressibility assumption, these inclusions give rise 
to two embeddings of $C$ into $A$, and the fundamental group of $X$
is the associated $HNN$ extension $A\ast _{C}$. 

\bigskip
\begin{theorem}
Let $\Gamma$  be the fundamental group of a compact riemannian manifold $X$
of dimension $n$ and sectional curvature less than or equal to $-1$.
Suppose that $\Gamma = A\ast_{C}$ where $A$ is a proper subgroup of $\Gamma$.
Then, we have $\delta(C) \geq n-2$ and equality $\delta(C) = n-2$ if and only if
there exist a non separating compact totally geodesic hypersurface $Y\subset X$
with fundamental group $\pi_{1}(Y) = C$. Moreover, in that case,
$Y$ with its induced metric is of constant sectional curvature $-1$, and the 
$HNN$ decomposition arising from $Y$ is the one we started with.  
\end{theorem}

Let us summarize the ideas of the proof of theorem 1.2.
We work on $\tilde{X}/C$. 
The amalgamation assumption provides an essential hypersurface $Z$ in 
$\tilde{X}/C$, namely $Z$ is homologically non trivial in $\tilde{X}/C$.
The volume of all hypersurfaces homologous to $Z$ is bounded below
by a positive constant because their systole are bounded away from zero.
We then construct a smooth map $F: \tilde{X}/C \rightarrow \tilde{X}/C$,
homotopic to the identity  
which contracts the volume of all compact 
hypersurfaces $Y$ by the factor $\big(\frac{\delta(C)}{n-2}\big)^{n-1}$,
namely $vol_{n-1}F(Y) \leq \big(\frac{\delta(C)}{n-2}\big)^{n-1}vol_{n-1} Y$.
This contracting property together with the lower bound of the volume
of hypersurfaces in the homology class of $Z$ gives
the inequality $\delta(C) \geq n-2$.
This map is different from the map constructed in \cite{bcg3}, in particular
it can be defined under the single condition that the limit set of $C$ is not reduced to one point.
Moreover, its derivative has an upper bound depending only on the critical exponent of $C$.

The equality case goes as follows.
When $\delta(C)=n-2$, the map $F: \tilde{X}/C \rightarrow \tilde{X}/C$ contracts 
the $(n-1)$-dimensional volumes, ie. $|Jac_{n-1}F|\leq 1$.
This contracting property is infinitesimally rigid in the following sense.
Let us consider a lift $\tilde{F}$ of $F$. If $|Jac_{n-1}\tilde{F}(x)| = 1$
at some point $x\in \tilde{X}$, then $\tilde{F}(x) = x$, there exists
a tangent hyperplane $E \subset T_{x}\tilde{X}$ such that 
$D\tilde{F}(x)$ is the orthogonal projector of $T_{x}\tilde{X}$ onto $E$
and the limit set $\Lambda_{C}$ is contained in the {\it topological equator}
$E(\infty)\subset \partial \tilde{X}$ associated to $E$. By topological
equator $E(\infty)\subset \partial \tilde{X}$ associated to $E$,
we mean the set of end points of those geodesic rays starting at $x$ tangently
to $E$. 

We then prove the existence of a point $x\in \tilde{X}$ such that 
\begin{equation}
|Jac_{n-1}\tilde{F}(x)| = 1.
\end{equation}
 If there would exist a minimizing 
cycle in the homology class of $Z$ in $\tilde{X}/C$, any point of such a cycle
would satisfy (1.1). As no such minimizing cycle a priori exists because of
non compactness of $\tilde{X}/C$, we prove instead the existence of a 
$L^{2}$ harmonic $(n-1)$-form dual to $Z$, which is enough to prove 
existence of a point $x$ such that (1.1) holds.
  
At this stage of the proof, 
there is a big difference between the constant curvature case
 and the variable curvature
case. 

In the constant curvature case, any topological equator bounds a totally geodesic
hyperbolic hypersurface $\Bbb {H}^{n-1}$, and therefore, as the group $C$ preserves
$\Lambda_{C} \subset E(\infty) = \partial\Bbb {H}^{n-1}$,
it is not hard to see that $C$ also preserves $\Bbb {H}^{n-1}$ and acts cocompactly on it,
and the hypersurface of the equality case in theorem 1.2 is $\Bbb {H}^{n-1}/C$, \cite{bcg2}.

In the variable curvature case, we first show the existence of a $C$-invariant
 totally geodesic hypersurface
$\tilde{Z}_{\infty} \subset \tilde{X}$ whose boundary at infinity coincides with $\Lambda(C)$, and
then we show that $\tilde{Z}_{\infty}$ is isometric to the real hyperbolic space.
We then show that $Y=:\mathbb{H}^{n-1}_{\mathbb{R}} /C$, which is compact, injects
in $X =\tilde{X}/\Gamma$ and separates $X$ in two connected components whose fundamental
groups are $A$ and $B$ respectively.

In order to show the existence of such a totally geodesic hypersurface $\tilde{Z}_{\infty}$,
we first prove that $C$ is a convex cocompact group, ie. the convex hull of the limit set 
of $C$ in $\tilde{X}$ has a compact quotient under the acton of $C$, and that 
the limit set of $C$ is homeomorphic to an $(n-2)$-dimensional sphere.

The convex cocompactness property of $C$ and the fact that 
the limit set $\Lambda(C)$ of $C$ is homeomorphic to an $(n-2)$-dimensional topological
sphere are the two key points in the equality case. 

 This compactness property then
allows us to prove the existence of a mimimizing current in the homology class of the essential
hypersurface $Z\subset \tilde{X}/C$. By regularity theorem this minimizing current $\tilde{Z}_{\infty}$ 
is a smooth manifold 
except at a singular set of codimension at least 8. By the contracting properties of our map $F$,
 $\tilde{Z}_{\infty}$ is fixed by $F$ and the geometric properties of $F$ at fixed points
where the $(n-1)$-jacobian of $F$ equals 1 allows us to prove that 
$\tilde{Z}_{\infty}$ is totally geodesic and isometric to the hyperbolic space.

Let us now briefly describe the proof of the convex cocompactness
property of $C$ in the equality case.

The group $C$ (or a finite index subgroup of it) actually globally preserves a smooth cocompact
hypersurface $\tilde{Z}\subset \tilde{X}$ which separates 
$\tilde{X}$ into two connected components and whose boundary 
$\partial \tilde{Z} \subset \partial \tilde{X}$
coincides with  $\Lambda_{C}\subset E(\infty)$.
In the case where  $C$ wouldn't be convex cocompact,
we are able to find an horoball $HB(\theta_{0})$ centered
at some point $\theta_{0} \in \Lambda_{C}$ in the complementary of which 
lies the hypersurface $\tilde{Z}$.

The contradiction then comes from the following.

Consider a sequence of points $\theta_{i} \in \partial \tilde{X}$
converging to $\theta_{0}$ and geodesic rays $\alpha_{i}$ starting from 
the point $x \in \tilde{X}$ at which $|Jac_{n-1}(x)| =1$ and ending up at
$\theta_{i}$.
These geodesic rays have to cross $\tilde{Z}$ at points $z_{i}$ which
are at bounded distance from the orbit $Cx$ of $x$, therefore the shadows
 $\mathcal{O}_{i}$
of balls centered at these $z_{i}$ enlighted from $x$ have to contain points 
of $\Lambda_{C}$ by the shadow lemma of D. Sullivan.
On the other hand, we show that it is possible to choose 
the sequence $\theta_{i}$ in such a way that these shadows $\mathcal{O}_{i}$
don't meet $\Lambda_{C}$. This property 
$\mathcal{O}_{i} \cap \Lambda_{C} = \emptyset$ comes from a choice of 
$\theta_{i}$ such that the distance between $z_{i}$ and the set $H$ of
all geodesics rays at $x$ tangent to $E\subset T_{x}\tilde{X}$
tends to $\infty$. Intuitively, in order to chose $z_{i}$ as far 
as possible from $H$,
the points $\theta_{i}$ have
to be chosen tranversally to $\Lambda_{C}$.
This transversality condition is not well defined because the limit
set $\Lambda_{C}$ might be highly non regular.
Thus, in order to prove that such a choice is possible,
we argue again by contradiction. 
If for any choice of a sequence 
$\theta_{i}$ converging to $\theta_{0}$, 
the distance between $z_{i}$ and $H$
stays bounded, then the Gromov distances $d(\theta_{i},\theta_{0})$ 
between $\theta_{i}$ and $\theta_{0}$ satisfy 
 $d(\theta_{i},\Lambda_{C})=o(d(\theta_{i},\theta_{0}))$,  
and therefore any tangent cone of $\Lambda_{C}$ at $\theta_{0}$
would coincide with a tangent cone of $\partial \tilde{X}$ at $\theta_{0}$,
which is known to be topologically $\Bbb{R}^{n-1}$.
But on the other hand, the existence of a point $x$ such that 
$|Jac_{n-1}(x)| =1$ and the fact that $C$ acts uniformly quasiconformally
with respect to the Gromov distance on $\partial\tilde{X}$ imply
that the Alexandroff compactification of the above tangent cone 
of $\Lambda_{C}$ at $\theta_{0}$ is homeomorphic to $\Lambda_{C}$
which is contained in a topological sphere $S^{n-2}$,
 leading to a contradiction.

From convex cocompactness of $C$ and the fact that the limit set of $C$ is a topological
$(n-2)$-dimensional sphere, there is an alternative proof of the existence
of a totally geodesic $C$-invariant copy of the hyperbolic space
 $\mathbb{H}^{n-1}_{\mathbb{R}} \subset \tilde{X}$
which consists in observing that the topological dimension and the 
Hausdorff dimension of the limit
set $\Lambda(C)$ are equal to $n-2$ 
and then use the following result of 
M. Bonk and B. Kleiner (which we quote in the riemannian manifold setting although it
remains true for $CAT(-1)$ spaces) instead of the (simpler) minimal current argument.

\begin{theorem} \cite{bk}
Let $X$ be a Cartan Hadamard 
$n$-dimensional manifold whose sectional curvature satisfy
$K \leq -1$, and $C$ a convex cocompact discrete subgroup
of isometries of $X$ with limit set $\Lambda_{C}$.
 Let us assume that the topological dimension 
and the Hausdorff dimension  
(with respect to the Gromov distance on 
$\partial \tilde{X}$) of $\Lambda_{C}$ coincide and are equal
 to an integer $p$.
 Then, $C$ preseves a totally geodesic embedded copy
of the real hyperbolic space $\Bbb{H}^{p+1}$, with
$\partial \mathbb H^{p+1} = \Lambda_{C}$. 
\end{theorem}

\bigskip
The authors would like to express their gratitude to Alex Lubotsky, Jean Barge, Marc Bourdon, Gilles Carron, Jean Lanne,
Fr\'{e}d\'{e}ric Paulin, Leonid Potyagailo for their interests and helpfull conversations. 

\section{Essential hypersurfaces}
Let $\Gamma$ be a discrete cocompact group of isometries of 
a $n$-dimensional Cartan-Hadamard manifold $(\tilde{X},\tilde{g})$
whose sectional curvature satisfies
$K_{\tilde{g}} \leq -1$. 
Let us assume that the compact manifold $X= \tilde{X}/\Gamma$ is orientable.
Let us also assume that 
$\Gamma = A \ast_{C} B$ is an amalgamated product of 
its subgroups $A$ and $B$ over $C$.

We first reduce to the case where 
 $[\Gamma : C]$ is 
infinite.
 
Namely, if $[\Gamma : C] < \infty$, then the critical exponent
$\delta(C) = \delta (\Gamma) \geq n-1$, and the equality in theorem 1.2
holds.

We then can assume that  $[\Gamma : C] = \infty$. 

\begin{lemma}
Let $\Gamma = A \ast_{C} B$ be as above, with 
 $[\Gamma : C] = \infty$. If neither $A$ nor $B$ equals $\Gamma$, then
 $H_{n-1}(\tilde{X} /C, \mathbb {Z}) \neq 0$.
\end{lemma}

{\bf Proof :}
 The Mayer-Vietoris sequence coming from the 
decomposition $\Gamma = A \ast_{C} B$ writes
 cf. \cite{br}, Corollary 7.7,

$$
 H_{n}(\tilde{X}/C,\Bbb {Z})
\rightarrow
H_{n}(\tilde{X}/A,\Bbb {Z})\oplus  H_{n}(\tilde{X}/B,\Bbb {Z})
\rightarrow
 H_{n}(\tilde{X}/\Gamma,\Bbb {Z})
\rightarrow...
$$

$$
...\rightarrow  H_{n-1}(\tilde{X}/C,\Bbb {Z})
\rightarrow ...
$$

$$
$$

As  $[\Gamma : C] = \infty$,
$ H_{n}(\tilde{X}/C,\Bbb {Z})=0$ thus, if 
 $H_{n-1}(\tilde{X} /C, \mathbb {Z}) = 0$,
 we deduce from
 the Mayer-Vietoris sequence that
$H_{n}(\tilde{X}/A,\Bbb {Z})\oplus  H_{n}(\tilde{X}/B,\Bbb {Z})$ is 
isomorphic to
$ H_{n}(\tilde{X}/\Gamma,\Bbb {Z})$. As
  $H_{n-1}(\tilde{X} /\Gamma, \mathbb {Z}) = \mathbb {Z}$,
we then deduce that either
$[\Gamma: A] = \infty$ and $B=\Gamma$, or $[\Gamma: B]=\infty$ and
$A=\Gamma$. $\Box$

\bigskip
In fact in the sequel of the paper
we will make use of a smooth essential hypersurface $Z$ in 
$\tilde{X}/C$.

\begin{definition}
A compact smooth orientable hypersurface $Z$ of
 an $n$-dimensional manifold $Y$  
is essential in $Y$ if
$i_{\ast}([Z]) \not= 0$
where $[Z]\in H_{n-1}(Z,\mathbb{R})$ denotes the fundamental class of $Z$
and 
$i_{\ast} :  H_{n-1}(Z,\Bbb {R}) \rightarrow  H_{n-1}(Y,\Bbb {R})$  
the morphism induced by the inclusion 
$i:Z \hookrightarrow Y$. 
\end{definition}

The end of this section is devoted to finding such an hypersurface $Z$
in $\tilde{X}/C$.

 Let us recall a few facts about amalgamated products and their actions
on trees, following \cite{se}. Let  
 $\Gamma = A \ast_{C} B$ be an amalgamated products of its subgroups
 $A$ and $B$ over $C$. Then, $\Gamma$ 
acts on a simplicial tree 
$\tilde{T}$ with a fundamental domain $T \subset \tilde{T}$ being 
a segment, ie. an edge joining two vertices. 
Let us describe this tree $\tilde{T}$.
There are two orbits of vertices $\Gamma v_{A}$ and $\Gamma v_{B}$,
the stabilizer of the edge $v_{A}$
(resp. $v_{B}$) being $A$,( resp. $B$).
There is one orbit of edges $\Gamma e_{C}$, the stabilizer of 
the edge $e_{C}$ being $C$.
The fundamental domain $T$ can be chosen as the edge $e_{C}$
joining the two vertices $v_{A}$ and $v_{B}$.
The set of vertices adjacent to $v_{A}$,
(resp. $v_{B}$), is in one to one correspondance with $A/C$, (resp. $B/C$).
Note that as neither $A$ nor $B$ are equal to $\Gamma$, then
 $[A : C]  \neq 1$ and $[B : C]  \neq 1$, therefore 
for an arbitrary point $t_{0}$ on the edge $e_{C}$ we see that 
$\tilde{T} - t_{0}$ is a disjoint union of two unbounded connected components.
This fact will be used later on. 

\bigskip
Let us consider a continuous $\Gamma$-equivariant map 
$\tilde{f}: \tilde{X} \rightarrow T$
where $T$ is the Bass-Serre tree associated to the amalgamation
$\Gamma = A \ast_{C} B$.
One regularizes $\tilde{f}$ such that it is smooth in restriction to the 
complementary of the inverse image of the 
set of vertices of $T$.
Let $t_{0}$ a regular value of  $\tilde{f}$ contained in 
that edge of $T$ which is fixed by the subgroup $C$ and define 
$\tilde{Z} = \tilde{f}^{-1}(t_{0})$.
$\tilde{Z}$ is a smooth orientable possibly not connected 
hypersurface in $\tilde{X}$, globally 
$C$-invariant.
Let us write $Z = \tilde{Z}/C$. 
We will show  $Z \subset \tilde{X}/C$ is compact and that one of
 the connected components of $Z$ is essential.

\begin{lemma}
 $Z \subset \tilde{X}/C$ is compact.
\end{lemma}  

{\bf Proof :} Let us show that for any sequence $z_{n} \in \tilde{Z}$,
there exists a subsequence $z_{n_{k}}$ and $\gamma _{k} \in C$ 
such that $\gamma _{k} z_{n_{k}}$ converges.
As $\Gamma$ is cocompact, there exists $g_{n}\in \Gamma$ such that 
the set $(g_{n} z_{n})$ is relatively compact. Let 
$g_{n_{k}} z_{n_{k}}$ a subsequence which converges to a point
 $z\in \tilde{X}$.
By continuity, the sequence $\tilde{f}(g_{n_{k}} z_{n_{k}})$ converges to 
$\tilde{f}(z)$, and by equivariance we get

$$
g_{n_{k}}\tilde{f}( z_{n_{k}})=g_{n_{k}} t_{0} \rightarrow \tilde{f}(z)
$$
when $k$ tends to $\infty$.
As $\Gamma$ acts in a simplicial way on the tree $T$ and transitively on the
set of edges,
the sequence 
$g_{n_{k}} t_{0}$ is stationary, ie 
$g_{n_{k}} t_{0} = t'_{0} = g t_{0}$ for $k$ large enough.
Thus $g^{-1}g_{n_{k}} = \gamma_{k} \in C$ for $k$ large enough
since it fixes $t_{0}$
and $\gamma_{k}  z_{n_{k}} = g^{-1}g_{n_{k}}  z_{n_{k}}$
converges to $g^{-1}(z)$ $\Box$

The smooth compact hypersurface $Z$ we constructed might be
not connected. Let us write 
$Z= Z_{1} \cup Z_{2} \cup... \cup Z_{k}$ where the 
$Z_{j}$'s are the connected components of $Z$. Each 
$Z_{j}$ is a compact smooth oriented hypersurface
of $\tilde{X}/C$.

The aim of what follows is to prove that at least one 
component $Z_{i}$ of $Z$ is essential.

\begin{lemma}
There exists $i\in [1,k]$ such that 
$Z_{i}$ is essential in $\tilde{X}/C$. 
\end{lemma}

{\bf Proof :} If there exists a $Z_{i}$ which doesn't separate 
$\tilde{X}/C$ in two connected components, then  $Z_{i}$ is essential
in $\tilde{X}/C$.
 So we can assume that every  $Z_{j}$, $j=1,...k$, does separate
$\tilde{X}/C$ in two connected components. In that case we will show
that there exists a $Z_{i}$ which separates 
$\tilde{X}/C$ in two unbounded connected components which easily 
implies that $Z_{i}$ is essential.

Let us denote $U_{l}$, $l=1,2,...,p$, the connected components of
$\tilde{X}/C - \cup_{j=1}^{k} Z_{j}$.

Claim : at least two components $U_{m}$, $U_{m'}$ are unbounded.

Assuming the claim let us finish the proof of the lemma.
For each $Z_{j}$ we denote $V_{j}$,  $V'_{j}$ the two connected components of
$\tilde{X}/C - Z_{j}$. Then 
$U_{m} = W_{1} \cap W_{2} \cap ...\cap W_{k}$ where for each $j$,
$W_{j} = V_{j}$ or $W_{j} = V'_{j}$. 
In the same way, 
$U_{m'} = W'_{1} \cap W'_{2} \cap ...\cap W'_{k}$.
As $U_{m}\cap U_{m'} = \emptyset$, there exists 
$i \in [1,k]$ such that $W_{i} \cap W'_{i} = \emptyset$, thus 
$U_{m}\subset V_{i}$ and  $U_{m'}\subset V'_{i}$ or 
$U_{m}\subset V'_{i}$ and  $U_{m'}\subset V_{i}$ so $Z_{i}$
separates $\tilde{X}/C$ into two unbounded components. 
This proves the lemma.

Let us prove the claim.

We have already noticed that $T - \{t_{0}\}$ is the disjoint union of
two unbounded connected components $T_{1}$ and  $T_{2}$. 
As $C$ acts on $T$ isometrically and simplicially then
$T/C - \{t_{0}\}= T_{1}/C \cup T_{2}/C$ is the disjoint union
of two unbounded connected components. 
Let $\bar{f} : \tilde{X}/C \rightarrow T/C$
the quotient map of $\tilde{f}$.
For each component $U_{i}$, we have 

$\bar{f}(U_{i}) \subset T_{1}/C$ or
  $\bar{f}(U_{i}) \subset T_{2}/C$, thus we can conclude 
the claim because $\bar{f}$ is onto $\Box$

Let $\pi : \tilde{X} \to \tilde{X}/C$ be the natural projection.
For any $i =1,2,...,k$, let us denote 
$\{ \tilde{Z}_{i}^{j}\}_{j\in J}$ the set of connected components
of $ \tilde{Z}_{i} =: \pi^{-1}(Z_{i})$.

For each $i\in [1,k]$, we claim that $C$ acts transitively
on the set $\{ \tilde{Z}_{i}^{j}\}_{j\in J}$. Namely,
let us consider  
$ \tilde{Z}_{i}^{j}$
$ \tilde{Z}_{i}^{j'}$,
$\tilde{z} \in  \tilde{Z}_{i}^{j}$,
$\tilde{z'} \in \tilde{Z}_{i}^{j'}$, 
and write $z=\pi \tilde{z}\in Z_{i}$ and $z'=\pi \tilde{z}'\in Z_{i}$.
Let $\alpha$ be a continuous path on $Z_{i}$ such that 
$\alpha(0) = z$ and $\alpha(1)= z'$, and $\tilde{\alpha}$ 
the lift of $\alpha$ such that $\tilde{\alpha}(0) = \tilde{z}$.
We have $\pi \tilde{\alpha}(1) = z'$ and
 $\tilde{\alpha}(1) \in \tilde{Z}_{i}^{j}$ for some $j$,
thus, there exists $c\in C$ such that 
$c(\tilde{\alpha}(1)) = \tilde{z}'$
 and therefore 
$ c\tilde{Z}_{i}^{j} =\tilde{Z}_{i}^{j'}$.$\Box$

\bigskip
Let us denote $C_{i}^{j}$ the stabilizer of $\tilde{Z}_{i}^{j}$,
and
 $Z_{i}^{j}= \tilde{Z}_{i}^{j}/C_{i}^{j} \subset \tilde{X}/C_{i}^{j}$.
Let us write 
$p :  \tilde{X}/C_{i}^{j} \to  \tilde{X}/C$ 
the natural projection.

\begin{lemma}
The restriction of $p$ to $Z_{i}^{j}$ is a diffeomorphism onto $Z_{i}$.
In particular,  $Z_{i}^{j}$ is compact. 
\end{lemma}

{\bf proof :} Let $z$ and $z'$ be two points in $Z_{i}^{j}$ such that 
$p(z) = p(z')$. 
Let $\tilde{z}$ and $\tilde{z'}$ be lifts of $z$ and $z'$ in 
$\tilde{X}$.
These two points $\tilde{z}$ and $\tilde{z'}$  which are in $\tilde{Z}_{i}$
actually belong to the same connected component $\tilde{Z}_{i}^{j}$
because for $j \neq j'$, 
$\tilde{Z}_{i}^{j'}/C_{i}^{j} \cap \tilde{Z}_{i}^{j}/C_{i}^{j} = \emptyset$.
As $p(z) = p(z')$, there exits $c\in C$ such that 
$\tilde{z}' = c \tilde{z}$, thus $c \in C_{i}^{j}$, and
$z=z'$, therefore the restriction of $p$ to $Z_{i}^{j}$ is injective.
The surjectivity comes from the fact that 
$\pi ^{-1}Z_{i} = \bigcup _{j\in J} \tilde{Z}_{i}^{j}$ and 
$C$ acts transitively on the set $\{ \tilde{Z}_{i}^{j}\} _{j\in J}$. $\Box$

Let us consider the integer $i\in [1,k]$ as in lemma 2.4,
ie. such that $ Z_{i} \hookrightarrow \tilde{X} /C$ is essential,
and choose $\tilde{Z}_{i}^{l}$ one component of $\pi^{-1}(Z_{i})$.

After possibly replacing $C_{i}^{l}$ by an index two subgroup, we may assume
that  $C_{i}^{l}$ globally preserves each of the two connected 
components $U_{i}^{l}$ and $V_{i}^{l}$ of 
$\tilde{X} -  \tilde{Z}_{i}^{l}$.  

\bigskip
\begin{lemma}
Let $i$, $l$ and $C_{i}^{l}$ be chosen as above. The compact hypersurface 
$Z_{i}^{l}=\tilde{Z}_{i}^{l}/C_{i}^{l}$ is essential
in $\tilde{X}/C_{i}^{l}$.  Moreover
the two connected components 
$U_{i}^{l}/ C_{i}^{l}$ and $V_{i}^{l}/C_{i}^{l}$ of 
$\tilde{X}/C_{i}^{l} -Z_{i}^{l}$ are unbounded.  
\end{lemma}

{\bf Proof :} 
Let us consider $p: \tilde{X}/C_{i}^{l} \rightarrow \tilde{X}/C$.
By lemma 2.5, the restriction of $p$ to $Z_{i}^{l}$ is a diffeomorphism
onto $Z_{i}$, therefore $Z_{i}^{l}$ is essential in  $\tilde{X}/C_{i}^{l}$
because $Z_{i}$ is essential in $\tilde{X}/C$.
As $C_{i}^{l}$  preseves $U_{i}^{l}$ and $V_{i}^{l}$,
$Z_{i}^{l}$ separates 
$\tilde{X}/C_{i}^{l}$ into two connected components
$U_{i}^{l}/ C_{i}^{l}$ and $V_{i}^{l}/C_{i}^{l}$.
and as 
$Z_{i}^{l}$ is essential in $\tilde{X}/C_{i}^{l}$,
$U_{i}^{l}/ C_{i}^{l}$ and $V_{i}^{l}/C_{i}^{l}$ are unbounded. $\Box$

{\bf In the sequel of the paper we will denote 
$\tilde{Z'} = \tilde{Z}_{i}^{l}$, $C'= C_{i}^{l}$ and 
$Z' = Z_{i}^{l} =  \tilde{Z}_{i}^{l}/C_{i}^{l}$.}

\section{Isosystolic inequality}
In this section we summarize facts and results due to
M.Gromov, \cite{gr}.
Let $Z$ be a $p$-dimensional compact orientable manifold and 
$i : Z \hookrightarrow Y$
an embedding of $Z$ into $Y$ where $Y$ is an aspherical space.
We suppose that $i_{\ast}([Z]) \not= 0$ where 
$i_{\ast} :  H_{p}(Z,\Bbb {R}) \rightarrow 
 H_{p}(Y,\Bbb {R})$ is  
the morphism induced by the embedding 
$Z \hookrightarrow Y$.
Let us fix a riemanniann metric $g$ on $Z$.
For each $z\in Z$ we consider the set $\mathcal{C}_{z}$ of those 
loops $\alpha$ at $z$ such that $i \circ \alpha$ is homotopically non trivial
in $Y$.

Let us define the systole of $(Z, g, i)$ at the point $z$ by
\begin{definition}
$sys_{i}(Z, g, z) = inf\{ lengh (\alpha)$, $\alpha \in \mathcal{C}_{z} \}$
\end{definition}
\noindent 
and the systole of  $(Z, g, i)$ by
\begin{definition}
$sys_{i}(Z, g) = inf\{sys_{i}(Z, g, z)$ $z\in Z\}$.
\end{definition}  

The following isosystolic inequality, due to M.Gromov says that the volume 
of any essential submanifold 
$Z$ of an aspherical space $Y$ relatively to any riemanniann 
metric on $Z$ is universally bounded below by it's systole.

\begin{theorem} \cite{gr}
There exists a constant $C_{p}$
such that for each p-dimensional riemanniann manifold $(Z,g)$ and any 
embedding 
$Z \hookrightarrow Y$ into an aspherical space $Y$ such that
 $i_{\ast}([Z]) \not= 0$ where 
$i_{\ast} :  H_{p}(Z,\Bbb {R}) \rightarrow 
 H_{p}(Y,\Bbb {R})$ is  
the induced morphism in homology,
then $vol_{p}(Z,g) \geq C_{p}(sys_{i}(Z,g))^{p}$
\end{theorem}

\bigskip
We will apply this volume estimates to the essential hypersurface
$i : Z \hookrightarrow \tilde{X}/C$ that we constructed in lemma 2.6.

The following lemma is immediate.
\begin{lemma}
Let  $C$ be a discrete
group acting on a simply connected manifold $\tilde{X}$, $\tilde{Z}$ a 
$C$-invariant hypersurface of $\tilde{X}$ and
$i : Z = \tilde{Z}/C \hookrightarrow \tilde{X}/C$
the natural inclusion.
Let $g$ any riemanniann metric
on $Z$ and $\tilde{g}$ the lift of $g$ to $\tilde{Z}$.
Then, for any $z\in Z$ we have,

 $$
sys_{i}(Z,g,z) = inf\{ d_{\tilde{g}}(\tilde{z},\gamma \tilde{z}), 
\gamma \in C \}
$$

\noindent
where $\tilde{z}\in \tilde{Z}$ is a lift of $z \in Z$ and 
$d_{\tilde{g}}$ is the distance induced by $\tilde{g}$ on $\tilde{Z}$.
\end{lemma}

{\bf Proof :} Let $\alpha \in \mathcal{C}_{z}$ a loop based at $z\in Z$.
As $i\circ \alpha$ is an homotopically non trivial loop at $i(z) = z$
in $\tilde{X}/C$, its lift $\widetilde{i\circ \alpha}$ at some
 $\tilde{z} \in \tilde{Z}$ ends up at $\gamma \tilde{z}$ for some
$\gamma \in C$.$\Box$

\section{Volume of hypersurfaces in $\tilde{X}/C$}
Let $(\tilde{X},\tilde{g})$ be a 
$n$-dimensional Cartan-Hadamard manifold 
whose sectional curvature satisfies $K_{\tilde{g}} \leq -1$ and
$C$ a discrete group of isometries of $(\tilde{X},\tilde{g})$.
We assume that the group $C$ is non elementary, namely $C$  fixes 
neither one nor two points in the geometric 
boundary $\partial \tilde{X}$ of $(\tilde{X},\tilde{g})$.

The aim of this section is to construct a map 
$F : \tilde{X}/C \rightarrow \tilde{X}/C$ 
such that for any compact hypersurface $Z$ of $\tilde{X}/C$, we
have 
$$
vol_{n-1}(F(Z)) \leq \Big(\frac{\delta +1}{n-1}\Big)^{n-1}vol_{n-1}(Z)
$$
where $\delta$ is the critical exponent of $C$ and $vol_{n-1}(Z)$ stands
for the $(n-1)$-dimensional volume of the metric on $Z$ induced from $g$.
For every subgroup $C' \subset C$ and any hypersurface 
$Z'$ of $\tilde{X}/C'$ 
the lift $F' : \tilde{X}/C' \rightarrow \tilde{X}/C'$ of $F$
will also verify 
$$
vol_{n-1}(F'(Z')) \leq \Big(\frac{\delta +1}{n-1}\Big)^{n-1}vol_{n-1}(Z').
$$

In order to construct the map
$F$ we need a few prelimiraries. We consider a finite positive Borel measure $\mu$ on the boundary 
$\partial \tilde{X}$ whose support contains at least two points.
Let us fix an origin $o \in \tilde{X}$ and denote $B(x,\theta)$
the Busemann function defined for each $x \in \tilde{X}$ and 
$\theta \in \partial \tilde{X}$ by
$$
B(x,\theta) = \lim _{t\to \infty} dist(x, c(t)) - t
$$
where $c(t)$ is the geodesic ray such that 
$c(0) = o $ and $c(+\infty) = \theta$.

Let $\mathcal D _\mu : \tilde{X} \to \Bbb{R}$
the function defined by 

\begin{equation}
\mathcal D _\mu (y) = \int _{\partial \tilde{X}} e^{B(y,\theta)} d\mu(\theta)
\end{equation}

A computation shows that 

\begin{equation}
Dd \mathcal D _\mu (y)
 = \int _{\partial \tilde{X}} (Dd B(y,\theta) + DB(y,\theta)\otimes DB(y,\theta))
e^{B(y,\theta)} d\mu(\theta).
\end{equation}

\noindent
When $K_{\tilde{g}} \leq -1$ the Rauch comparison theorem says that for every $y\in \tilde{X}$,
and $\theta \in \partial \tilde{X}$,

\begin{equation}
Dd B(y,\theta) + DB(y,\theta) \otimes DB(y,\theta) \geq \tilde{g}.
\end{equation}

We then get

\begin{equation}
Dd \mathcal D _\mu (y) \geq \mathcal D _\mu (y)  \tilde{g},
\end{equation}

thus $Dd \mathcal D _\mu (y)$ is positive definite and $\mathcal D _ \mu$
is strictly convex.

\begin{lemma} We have
$\lim_{y_{k} \to \partial \tilde{X}} \mathcal D _\mu (y) = + \infty.$
\end{lemma}

{\bf proof :} Let $y_{k} \in \tilde{X}$ a sequence such that 
\begin{equation}
\lim_{k \to \infty} y_{k} = \theta _{0} \in \partial \tilde{X}.
\end{equation}

As $supp(\mu) \cap (\partial \tilde{X} -\{\theta _{0}\}) \neq \emptyset$,
there exists a compact subset $K \subset \partial \tilde{X} -\{\theta _{0}\}$
such that $\mu (K) > 0$ thus,

\begin{equation}
\int _{\partial \tilde{X}} e^{B(y_{k},\theta)} d\mu \geq 
\int _{K} e^{B(y_{k},\theta)} d\mu  \to +\infty. 
\end{equation}
$\Box$

\begin{cor}
Let $\mu$ a finite borel measure on $\partial \tilde{X}$ 
whose support contains at least two points. The function 
$\mathcal{D}_{\mu}$ has a unique minimum.
This minimum will be denoted by $\mathcal{C}(\mu)$.
\end{cor}

Let us now consider some discrete subgroup $C \subset Isom (\tilde{X},\tilde{g})$.
Recall that a family of Patterson measures 
$(\mu _{x})_{x\in \tilde{X}}$ 
associated to $C$ is a set of positive finite measures $\mu _{x}$ on
$\partial \tilde{X}$, $x\in \tilde{X}$, such that the following holds for all 
$x\in \tilde{X}$, $\gamma \in C$,

\begin{equation}
\mu _{\gamma x} = \gamma _{*} \mu _{x}
\end{equation}

\begin{equation}
\mu _{x} = e^{-\delta B(x,\theta)} \mu_{o},
\end{equation}

where $o\in \tilde{X}$ is a fixed origin, $B$ the Busemann function associated to $o$ 
and $\delta$ the critical exponent of $C$.

We assume now that $supp (\mu _{o})$ contains at least two points and define the 
map $\tilde{F} : \tilde{X} \to  \tilde{X}$ for
$x\in \tilde{X}$ by

\begin{equation}
\tilde{F}(x) = \mathcal{C} (e^{-B(x,\theta)} \mu_{x}). 
\end{equation}

Here are a few notations.
For a subspace $E$ of $T_{x}\tilde{X}$, we will write 
$Jac_{E} \tilde{F}(x)$ the determinant of the matrix of the 
restriction of $D\tilde{F}(x)$ to $E$ with respect to orthonormal
bases of $E$ and $D\tilde{F}(x)E$. For an integer
$p$, we denote by $Jac_{p}\tilde{F}(x)$ the supremum
of $|Jac_{E} \tilde{F}(x)|$ as $E$ runs through the set of $p$-dimensional
subspaces of  $T_{x}\tilde{X}$.

\begin{lemma}
The map $\tilde{F}$ is smooth, homotopic to the Identity
 and verifies for all $x\in \tilde{X}$, $\gamma \in C$
and $p \in [2, n= dim(X)]$,

\bigskip
\noindent
(i) $\tilde{F}(\gamma x) = \gamma \tilde{F}(x)$

\bigskip
\noindent
(ii)  $|Jac _{p}\tilde{F}(x)| \leq  \Big(\frac {(\delta +1)}{p}\Big)^{p}$. 
\end{lemma}

{\bf Proof :} 

The map
 
$$
(x, y) \to 
\int _{\partial \tilde{X}} e^{B(y,\theta) - B(x,\theta)}d\mu _{x}(\theta) =
\int _{\partial \tilde{X}} e^{B(y,\theta) - (\delta + 1) B(x,\theta)}d\mu _{o}(\theta)
$$

is smooth because $y \to B(y,\theta)$ is smooth.

For all $x$ the map 
$y \to
\int _{\partial \tilde{X}} e^{B(y,\theta) - (\delta + 1) B(x,\theta)}d\mu _{o}(\theta)$
is strictly convex by (4.4) and tends to infinity when $y$ tend to $\partial \tilde{X}$ (cf.
lemma 4.1),
thus the unique minimum 
$\tilde{F}(x)$  is a smooth function. 
 The equivariance of $\tilde{F}$ comes from the cocycle relation
$B(\gamma y, \gamma \theta) - B(\gamma x, \gamma \theta) =
 B(y,\theta) - B(x,\theta)$.

For each $x\in \tilde{X}$ let $c _{x}$ be the geodesic in $\tilde{X}$ such that 
$c_{x}(0) = x$, $c_{x}(1) = \tilde{F}(x)$ and which is parametrized with constant speed.
The map $\tilde{F}_{t} : \tilde{X} \to \tilde{X}$
defined by $\tilde{F}_{t}(x) = c_{x}(t)$ is a $C$-equivariant homotopy between
$Id_{\tilde{X}}$ and $\tilde{F}$.

It remains to  prove (ii).

The point $\tilde{F}(x)$ is characterized by 

\begin{equation}
\int _{\partial \tilde{X}} DB(\tilde{F}(x),\theta) 
e^{B(\tilde{F}(x),\theta) - B(x,\theta)}d\mu _{x}(\theta) = 0.
\end{equation} 

In order to simplify the notations 
we will write $B_{(x,\theta)}$ instead of 
$B(x,\theta)$ 
and we will denote $\nu_{x}$ the measure
$e^{B(\tilde{F}(x),\theta) - B(x,\theta)}\mu _{x}$.
We will also write $D\tilde{F}(u)$ instead of $D\tilde{F}(x)(u)$.

The differential of $ \tilde{F}$ is characterized by the following: for 
$u\in T_{x}\tilde{X}$ and  $v\in T_{\tilde{F}(x)}\tilde{X}$, one has

$$
\int _{\partial \tilde{X}} [DdB_{(\tilde{F}(x),\theta)}( D\tilde{F}(u),v)
+ DB_{(\tilde{F}(x),\theta)}(v) DB_{(\tilde{F}(x),\theta)}( D\tilde{F}(u))]
d\nu_{x}(\theta)
$$

\begin{equation}
=(\delta + 1) 
\int _{\partial \tilde{X}} DB_{(\tilde{F}(x),\theta)}(v)
 DB_{(x,\theta)}(u)
d\nu _{x}(\theta).
\end{equation}
  
We define the quadratic forms $k$ and $h$
 for $v\in T_{\tilde{F}(x)}\tilde{X}$ by

\begin{equation}
k(v,v)=
\int _{\partial \tilde{X}} [DdB_{(\tilde{F}(x),\theta)}(v,v)+
 (DB_{(\tilde{F}(x),\theta)}(v))^{2}]
d\nu _{x}(\theta).
\end{equation}

and

\begin{equation}
h(v,v)=\int_{\partial\tilde{X}} DB_{(\tilde{F}(x),\theta)}(v)^{2}d\nu_{x}(\theta).
\end{equation}

The relation (4.11) writes, for $u\in T_{x}\tilde{X}$  and $v\in T_{\tilde{F}(x)}\tilde{X}$ :

\begin{equation}
k(D\tilde{F}(u),v) = (\delta +1)
\int _{\partial \tilde{X}} 
DB_{(\tilde{F}(x),\theta)}(v) DB_{(x,\theta)}(u)
d\nu _{x}(\theta).
\end{equation}

We defines the quadratic form $h'$ on $T_{x}\tilde{X}$ for $u\in T_{x}\tilde{X}$ by

\begin{equation}
h'(u,u)=
\int _{\partial \tilde{X}} 
DB_{(x,\theta)}(u)^{2}
d\nu _{x}(\theta),
\end{equation} 

and one derives from (4.14)

\begin{equation}
|k(D\tilde{F}(x)(u),v)| \leq 
(\delta +1) h(v,v)^{1/2} h'(u,u)^{1/2}.
\end{equation}

One now can estimate $Jac_{p}\tilde{F}(x)$.
Let $P \subset T_{x}\tilde{X}$, $dim P =p$.
If  $D\tilde{F}(P)$ has dimension lower than $p$, then there is nothing to be proven.
Let us assume that $dim D\tilde{F}(P)= p$.
Denote by the same letters $H'$ [resp. $H$ and $K$] the selfadjoint operators
 (with respect to $\tilde{g}$)
associated to the quadratic forms $h'$ [resp. $h$, $k$]
restricted to $P$ [resp. $D\tilde{F}(P)$].

Let $(v_{i})_{i=1}^{p}$ an orthonormal basis of $D\tilde{F}(P)$ 
which diagonalizes $H$ and 
$(u_{i})_{i=1}^{p}$ an orthonormal basis of $P$ such that
the matrix of $K\circ D\tilde{F}(x) : P \rightarrow D\tilde{F}(P)$ is
triangular. Then,

\begin{equation}
det K .|Jac_{P}\tilde{F}(x)| \leq (\delta +1)^{p}
(\Pi _{i=1}^{p} h(v_{i},v_{i})^{1/2}) 
(\Pi _{i=1}^{p} h'(u_{i},u_{i})^{1/2}) 
\end{equation}

thus,

\begin{equation}
det K .|Jac_{P}\tilde{F}(x)| \leq (\delta +1)^{p}
\Big(\frac{Trace H}{p}\Big)^{p/2} \Big(\frac{Trace H'}{p}\Big)^{p/2}.
\end{equation}

In these inequalities one can normalize the measures 
$$
\nu_{x}= e^{B(\tilde{F}(x),\theta) - B(x,\theta)}\mu _{x}
$$
such that their total mass equals one, which gives

\begin{equation}
trace H = \Sigma _{i=1}^{p} h(v_{i},v_{i}) \leq 1,
\end{equation}

the last inequality coming from the fact that for all $\theta \in \partial \tilde{X}$,
\begin{equation}
\Sigma _{i=1}^{p} DB(\tilde{F}(x),\theta)(v_{i})^{2}
\leq
||\nabla B(\tilde{F}(x),\theta)||^{2}=1
\end{equation}

and from the previous normalization.

Similarly, 

\begin{equation}
trace H' = \Sigma _{i=1}^{p} h'(u_{i},u_{i}) \leq 1.
\end{equation}

We then obtain with (4.18)

\begin{equation}
det K .|Jac_{P} \tilde{F}(x)| \leq \Big(\frac{\delta +1}{p}\Big)^{p}.
\end{equation}

Thanks to (4.3), we have $det K \geq 1$, so

\begin{equation}
|Jac_{P} \tilde{F}(x)| \leq \Big(\frac{\delta +1}{p}\Big)^{p}.  
\end{equation}

We get (ii) by taking the supremum in $P$.
$\Box$

As the map $\tilde{F} : \tilde{X} \to \tilde{X}$ is $C$-equivariant, 
then for every subgroup $C'\subset C$, $\tilde{F}$ gives rise to 
a map $F' : \tilde{X} /C' \to \tilde{X} /C'$ and so does the homotopy $\tilde{F}_{t}$ between 
$\tilde{F}$ and $Id_{\tilde{X}}$.

\begin{cor}
The map $F' : \tilde{X}/C' \to \tilde{X}/C'$ is homotopic to the Identity 
map and verifies for all
$x\in \tilde{X}/C'$ and $p\in [2, n=dim X]$ 
$$
|Jac F'_{p}(x)| \leq \Big(\frac{\delta +1}{p}\Big)^{p} . 
$$  
\end{cor}

Let $C\subset Isom (\tilde{X},\tilde{g})$ as above, ie such that the support of the Patterson-
Sullivan measures contains at least two points and with critical exponent $\delta$.
Let $C' \subset C$ be a subgroup.

Let us consider an compact hypersurface $Z' \subset \tilde{X}/C'$.

Denote $F'^{k} = F'\circ F'\circ .....\circ F'$ the composition of $F'$ $k$-times.

Let us write $g_{k} = (F'^{k})^{*}g$, where $g$ is the metric on $\tilde{X}/C'$ 
induced by $\tilde{g}$.
The symmetric 2-tensor $g_{k}$ may not be a riemanniann metric on $\tilde{X}/C'$ nor 
its resriction 
to $Z'$, so we have to modify it. 
For $\epsilon > 0$, the following symmetric 2-tensor
$ g_{\epsilon,k}$ is a riemanniann 
metric on $\tilde{X}/C'$ and so is its restriction $h_{\epsilon,k}$ to $Z'$.

\begin{equation}
g_{\epsilon,k} = g_{k} + \epsilon ^{2}g.
\end{equation}

\begin{lemma}
Let $h_{\epsilon,k}$ be the restriction of  $g_{\epsilon,k}$ to the hypersurface $Z'$
and $g_{Z'}$ the restriction of $g$ to $Z'$.
Let  $\Phi_{\epsilon,k} : Z' \to \Bbb R$ the density defined for all $x\in Z'$ by
 $dv_{h_{\epsilon,k}}(x) = \Phi_{\epsilon,k} (x) dv_{g_{Z'}}(x)$.
For any sequence $\epsilon_{k}$ such that $\lim_{k \to \infty} \epsilon_{k} = 0$,
there exists a sequence $\epsilon_{k}'$, $\lim_{k \to \infty} \epsilon_{k}' = 0$,
such that for all $x\in Z'$,
$$
0 < \Phi_{\epsilon_{k}',k}(x) \leq |Jac _{n-1} F'^{k}(x)| + \epsilon_{k}.
$$
In particular, 
$$
 \Phi_{\epsilon_{k}',k}(x) 
\leq \Big(\frac{\delta +1}{n-1}\Big)^{k(n-1)}+ \epsilon_{k}
$$
and
$$
vol(Z',h_{\epsilon_{k}',k}) 
\leq \Big[\Big(\frac{\delta +1}{n-1}\Big)^{k(n-1)}+
 \epsilon_{k}\Big]vol(Z',g_{Z'}).
$$
\end{lemma}

\begin{cor}
Under the above asumptions, if $\delta < n-2$ there exists a sequence 
$\epsilon_{k}'$ such that $\lim_{k \to \infty} \epsilon_{k}' = 0$,
and
$ lim_{k\to \infty} vol(Z',h_{\epsilon_{k}',k}) = 0$.
\end{cor}

{\bf Proof of lemma 4.5 :}

Let us fix $k$ an integer.
Let $x\in Z'$ and $u\in T_{x}Z'$.
We have 
$g_{\epsilon,k}(u,u) = h_{\epsilon,k}(u,u) = 
g(DF'^{k}(x)(u),DF'^{k}(x)(u)) + \epsilon^{2}g(u,u)$
thus 
$h_{\epsilon,k}(u,u) = g(A_{x,\epsilon}u,u)$
where $A_{x,\epsilon} \in End (T_{x}Z')$ is the self adjoint operator
$A_{x} = DF'^{k}(x)^{*} \circ DF'^{k}(x) + \epsilon ^{2} Id$, 
with $ DF'^{k}(x)^{*}$ the adjoint of 
$ DF'^{k}(x)  : (T_{x}Z',g(x)) \to (DF'^{k}(x)(T_{x}Z'), g(F'^{k}(x))).$

By compactness of $Z'$ and continuity of $A_{x,\epsilon}$,
there exist $\epsilon'_{k}$ such that
$$
\Phi_{\epsilon'_{k} ,k}(x) = det A_{x,\epsilon'_{k}}^{1/2}
\leq det A_{x,0} + \epsilon_{k},
$$

thus
$$
\Phi_{\epsilon_{k}' ,k}(x) \leq |Jac _{n-1} F'^{k}(x)| + \epsilon_{k}. 
$$

The lemma then follows from Corollary 4.4. $\Box$

\section{Proof of Theorem 1.2} 
This section is devoted to the proof of the Theorem 1.2.
Let $\Gamma$ be a discrete cocompact group of isometries of 
a $n$-dimensional Cartan-Hadamard manifold $(\tilde{X},\tilde{g})$
whose sectional curvature satisfies
$K_{\tilde{g}} \leq -1$. 
We assume that 
$\Gamma = A \ast_{C} B$.
At the end of section 2 we constructed a subgroup $C' \subset C$ and an orientable hypersurface 
$\tilde{Z'}\subset \tilde{X}$
such that $C'. \tilde{Z'} = \tilde{Z'}$ and 
$Z' = \tilde{Z'}/C'$ is compact in $\tilde{X}/C'$.
Moreover $Z'$ is essential in $\tilde{X}/C'$ ie 
$i_{*}([Z']) \neq 0$ 
where 
$i_{*} : H_{n-1}(Z',\Bbb R) \to H_{n-1}(\tilde{X}/C',\Bbb R)$
is the morphism induced on homology groups by the inclusion
$i : Z' \to \tilde{X}/C'$
and $[Z']$ the fundamental class of $Z'$.

{\bf 5.1 Proof of the inequality}    

We now prove the inequality in the theorem 1.2.
Let us assume that $\delta < n-2$ and derive a contradiction.
Let $h_{\epsilon_{k}',k}$ the sequence of metric defined on $Z'$ in lemma 4.5,
then by corollary 4.6 we have 
\begin{equation}
\lim_{k\to \infty} vol(Z',h_{\epsilon_{k}',k}) = 0.
\end{equation}

We now show that the systole of the metric  $h_{\epsilon_{k}',k}$ on $Z'$ is bounded 
below independently of $k$.
Recall that the systole of $i : Z' \rightarrow \tilde{X}/C'$
at a point $z\in Z'$
with respect to a metric $h_{\epsilon , k}$ can be defined by

\begin{equation}
\rm{sys}_{i}(Z' , h_{\epsilon ,k}, z) =
inf_{\gamma \in C'}dist_{(\tilde{Z'},\tilde{h}_{\epsilon , k})}(\tilde{z},
 \gamma \tilde{z})
\end{equation}

where $\tilde{z}$ is any lift of $z$ and
$\tilde{h}_{\epsilon , k}$ the lift on $\tilde{Z'}$ of $h_{\epsilon , k}$,
(cf. lemma 3.4).
 
Let $\alpha(t)$ be a minimizing geodesic between $\tilde{z}$
 and $\gamma \tilde{z}$
on $(Z',\tilde{h}_{\epsilon , k})$.
By definition of $\tilde{h}_{\epsilon , k}$ we have

\begin{equation}
\rm{dist}_{(\tilde{Z'},\tilde{h}_{\epsilon , k})}(\tilde{z}, \gamma \tilde{z})
\geq l_{\tilde{g}}(\tilde{F}^{k}\circ \alpha)
\end{equation}

where
$l_{\tilde{g}}$
 stands for the lengh with respect to
 $\tilde{g}$
on 
$\tilde{X}$.

We get 

\begin{equation}
\rm{dist}_{(\tilde{Z'},\tilde{h}_{\epsilon , k})}(\tilde{z}, \gamma \tilde{z})
\geq \rm{dist}_{(\tilde{X},\tilde{g})}(\tilde{F}^{k}(\tilde{z}),\gamma \tilde{F}^{k}(\tilde{z}))
\geq \rho
\end{equation}

where $\rho$
is the injectivity radius of $\tilde{X}/C'$.

We then have

\begin{equation}
\rm{sys}_{i}(Z' , h_{\epsilon ,k}) \geq \rho,
\end{equation}

and thanks to the Theorem 3.3 we obtain

\begin{equation}
\rm{vol}(Z',h_{\epsilon'_{k},k}) \geq C_{n}\rho ^{n-1}
\end{equation}

which contradicts (5.1).

$\Box$

{\bf 5.2 Proof of the equality case}

There will be several steps.

\bigskip
\noindent
{\bf Step 1:} The limit set of $C$ is contained in a topological equator.

\bigskip
\noindent
{\bf Step 2:} The weak tangent to $\partial \tilde{X}$ and $\Lambda_{C'}$.

\bigskip
\noindent
{\bf Step 3:} The limit set $\Lambda_{C'}$ of $C'$
and the limit set $\Lambda_{C}$ of $C$ are equal to a topological equator.

\bigskip
\noindent
{\bf Step 4:} $C'$ and $C$ are convex cocompact.

\bigskip
\noindent
{\bf Step 5:} $C$ preserves a copy of the real hyperbolic space 
$\mathbb{H}^{n-1}_{\mathbb{R}}$ totally geodesically embedded in $\tilde{X}$.

\bigskip
\noindent
{\bf Step 6:} Conclusion

\bigskip
\noindent
{\bf Step 1: The limit set $\Lambda _{C}$ of $C$ is contained in a topological equator.} 

\bigskip
\noindent
Let $x\in \tilde{X}$
and $E\subset T_{x}\tilde{X}$
a codimension one subspace.
For each $u\in  T_{x}\tilde{X}$, $\tilde{g}(u,u)=1$, one considers the geodesic 
$c_{u}$ defined by 
$c_{u}(0) = x$
and $\dot{c_{u}}(0) = u$.
We define the equator $E(\infty)$ associated to $E$ as the subset of
$\partial \tilde{X}$

\begin{equation}
E(\infty) = \{ c_{u}(+\infty) / u\in E\}
\end{equation}

Our goal is to prove the existence of a point $x\in \tilde{X}$ and an 
hyperplane 
$E\subset T_{x}\tilde{X}$ such that the limit set $\Lambda _{C}$
satisfies  $\Lambda _{C} \subset E(\infty)$.

Recall that $C'\subset C$ globally preserves an hypersurface $\tilde{Z'}$
such that 
$\tilde{Z'}/C' \subset \tilde{X}/C'$
is  compact and essential.

Let us also recall that we have constructed a $C$-equivariant map
$\tilde{F}: \tilde{X} \to \tilde{X}$ such that, for all $x\in \tilde{X}$,

\begin{equation}
|Jac_{n-1}\tilde{F}(x)| \leq \Big( \frac{\delta +1}{n-1}\Big)^{n-1} 
\end{equation} 

where the critical exponent $\delta$ of $C$ satisfies $\delta =n-2$, thus

\begin{equation}
|Jac_{n-1}\tilde{F}(x)| \leq 1.
\end{equation} 

The step 1 follows from the two following Propositions.

\begin{prop}
Let $x\in \tilde{X}$ such that $|Jac_{n-1}\tilde{F}(x)| = 1$.
Then there exists $E\subset T_{x}\tilde{X}$ such that the limit set
 $\Lambda _{C}$
satisfies  $\Lambda _{C} \subset E(\infty)$.
Moreover, $\tilde{F}(x) =x$ and $D\tilde{F}(x)$ is the orthogonal
projector onto $E$.
\end{prop}

\begin{prop}
There exists $x\in \tilde{X}$ such that $|Jac_{n-1}\tilde{F}(x)| = 1$.
\end{prop}

\bigskip
\noindent
{\bf Proof of Proposition 5.1}

\bigskip
\noindent
Let $x\in \tilde{X}$
such that $|Jac_{n-1}\tilde{F}(x)| = 1$.
By definition we have  a subspace $E\subset T_{x}\tilde{X}$ such that 
$|Jac_{E}\tilde{F}(x)| = 1$.
By (4.18) and $\rm det K \geq 1$ we have,

$$
|Jac_{E}\tilde{F}(x)| \leq (n-1)^{n-1}
\Big(\frac{trace H}{n-1}\Big)^{\frac{n-1}{2}}
\Big(\frac{trace H'}{n-1}\Big)^{\frac{n-1}{2}}
$$

\begin{equation}
\leq (n-1)^{n-1}(\frac{1}{n-1})^{n-1}.
\end{equation}

In particular as $|Jac_{E}\tilde{F}(x)| = 1$, we have equality in the inequalities (5.10),
thus, $trace H = trace (h) =1$, and

\begin{equation}
H = \frac{1}{n-1} Id_{D\tilde{F}(x)(E)}.
\end{equation}

Let us recall that the quadratic form $h$ is defined by
$$
h(v,v)=
\int_{\partial\tilde{X}} DB(\tilde{F}(x),\theta)(v)^{2}
e^{B(\tilde{F}(x),\theta)-B(x,\theta)}d\mu_{x}(\theta)
$$
where
$\mu_{x}$ is the Patterson-Sullivan measure of $C$ 
normalized by

\begin{equation}
\int_{\partial\tilde{X}}
e^{B(\tilde{F}(x),\theta)-B(x,\theta)}d\mu_{x}(\theta)=1.
\end{equation}

We then have
$$
1= trace (h) = trace H =\Sigma_{i=1}^{n-1}h(v_{i},v_{i}) =
$$
$$
= \int_{\partial \tilde{X}}\Sigma_{i=1}^{n-1}DB(\tilde{F}(x),\theta)(v_{i})^{2}
e^{B(\tilde{F}(x),\theta)-B(x,\theta)}d\mu_{x}(\theta)
$$
$$
\leq  \int_{\partial \tilde{X}}||DB(\tilde{F}(x),\theta)||^{2} 
e^{B(\tilde{F}(x),\theta)-B(x,\theta)}d\mu_{x}(\theta) \leq 1,
$$

because 
$\Sigma_{i=1}^{n-1}DB(\tilde{F}(x),\theta)(v_{i})^{2} \leq
||DB(\tilde{F}(x),\theta)||^{2} =1$

for all $\theta \in \partial \tilde{X}$.

Therefore for $\mu_{x}$-almost all 
$\theta \in supp(e^{B(\tilde{F}(x),\theta)-B(x,\theta)}d\mu_{x}(\theta))
 = supp(\mu_{x})$, we have

\begin{equation}
\Sigma_{i=1}^{n-1}DB(\tilde{F}(x),\theta)(v_{i})^{2} =
||DB(\tilde{F}(x),\theta)||^{2} =1.
\end{equation}

In (5.12), $\Sigma_{i=1}^{n-1} DB(\tilde{F}(x),\theta)(v_{i})^{2}$
represents the square of the norm of the projection of $\nabla B(\tilde{F}(x),\theta)$ on $E$.

By continuity of $B(x,\theta)$ in $\theta$ one then gets 
$\Lambda_{C} = supp(\mu_{x}) \subset E(\infty)$.

Let us now prove that $\tilde{F}(x) =x$.
When $Jac \tilde{F}_{E}(x) =1$, we have equality 
in the Cauchy-Schwarz inequality (4.16), therefore
for each $i = 1,...n-1$ and $\theta \in \Lambda_{C}$ we get 
$DB(\tilde{F}(x),\theta)(v_{i}) = DB(x,\theta)(u_{i})$.
Therefore we deduces from (5.13) that  
$\nabla B(x,\theta)= \Sigma _{i=1}^{i=n-1} DB(x , \theta)(u_{i}) u_{i}$,
which imply with (4.10) that 
\begin{equation}
\int _{\partial \tilde{X}}
DB(x,\theta)e^{B(\tilde{F}(x),\theta)-B(x,\theta)}
d\mu_{x}(\theta) = 0.
\end{equation}
On the other hand, as  
$\int _{\partial \tilde{X}}
DB(\tilde{F}(x),\theta)e^{B(\tilde{F}(x),\theta)-B(x,\theta)}d\mu_{x}(\theta)
=0$
 and
$H=\frac{1}{n-1}Id_{D\tilde{F}(x)(E)}$, the support of $\mu_{x}$ cannot be 
just a pair of points, therefore
the barycenter of the measure 
$e^{B(\tilde{F}(x),\theta) - B(x, \theta)}\mu_{x}$
defined in \cite{bcg2} is well defined and characterized as the point $z\in \tilde{X}$
 such that
$$
\int _{\partial \tilde{X}}
DB(z,\theta)e^{B(\tilde{F}(x),\theta)-B(x,\theta)}
d\mu_{x}(\theta) = 0,
$$
 thus (5.14) and (4.10) imply $x=\tilde{F}(x)$.
$\Box$

{\bf Proof of Proposition 5.2 :}

\bigskip
\noindent
If we knew that there exists a minimizing hypersurface $Z_{0}$
in the homology class of $Z$, then every points $x\in Z_{0}$ would 
verify  $|Jac_{n-1}\tilde{F}(x)| = 1$. 
We unfortunately don't know if there exists such a minimizing hypersurface
nor a minimizing current in the homology class of $Z$. Instead we will consider  
an $L^{2}(\tilde{X}/C')$ harmonic $(n-1)$-form dual to the homology class of $Z$.

We need the following lemmas in order to prove the existence of such a dual form.

Let $\lambda _{1}(\tilde{X}/C')$ be the bottom of the spectrum of the 
Laplacian on $(\tilde{X},\tilde{g})$, ie.

\begin{equation}
\lambda _{1}(\tilde{X}/C') = 
inf _{u\in C^{\infty}_{0}(\tilde{X}/C')}
\{ \frac{\int_{\tilde{X}/C'}|du|^{2}}{\int_{\tilde{X}/C'}u^{2}}\}.
\end{equation} 

\begin{lemma}
Let $C \subset Isom(\tilde{X},\tilde{g})$ a discrete group of isometries 
with critical exponent $\delta = n-2$ where 
$(\tilde{X},\tilde{g})$ is an n-dimensional 
Cartan-Hadamard manifold of sectional curvature
$K_{\tilde{g}} \leq -1$.
Then for any subgroup $C' \subset C$ we have 
$\lambda _{1}(\tilde{X}/C') \geq n-2$.
\end{lemma}

{\bf Proof :}  

\bigskip
\noindent
Thanks to a theorem of Barta, cf.\cite{su} Theorem 2.1, the lemma boils down to finding  a 
positive function
$c : \tilde{X}/C' \to \Bbb R_{+}$ such that 
$\Delta c(x) \geq  (n-2)c(x)$. 
Here, the laplacian $\Delta$ is the positive operator ie. $\Delta c = -trace Dd c$. 
We consider the smooth function $\tilde{c} : \tilde{X} \to \Bbb R_{+}$
defined by 
$\tilde{c}(x) = \mu_{x}(\partial \tilde{X})$ 
where $\{\mu_{x}\}_{x\in \tilde{X}}$
is a family of Patterson-Sullivan measure of $C$.
The function $\tilde{c}$ is $C$-equivariant therefore it defines
a map $c : \tilde{X}/C' \to \Bbb R_{+}$ 
 for any subgroup
$C'\subset C$.
Let us show 

\begin{equation}
\Delta c(x) \geq \delta (n-1-\delta) c(x) = n-2.
\end{equation}  

We have 

$$
\tilde{c}(x) = \int_{\partial \tilde{X}} e^{-\delta B(x,\theta)}d\mu_{o}(\theta)
$$

therefore

$$
\Delta \tilde{c}(x) =
\int_{\partial \tilde{X}}[-\delta \Delta B(x,\theta) - \delta ^{2}]
d\mu_{x}(\theta).
$$ 

The sectional curvature 
 $K_{\tilde{g}}$ of 
$(\tilde{X},\tilde{g})$
satisfies  
$K_{\tilde{g}} \leq -1$  we thus have
$-\Delta B(x,\theta) \geq n-1$ and as $\delta = n-2$
 we get 

$$ 
\Delta \tilde{c}(x) \geq
[\delta (n-1) - \delta ^{2}]\tilde{c}(x) = (n-2)\tilde{c}(x).
$$

$\Box$

The following lemma is due to G.Carron and E.Pedon, \cite{cp}.
For a complete riemannian manifold $Y$, we denote $H_{c}^{1}(Y,\mathbb R)$ the first cohomology
group generated by diferential forms with compact support.

\begin{lemma} [\cite{cp}, Lemme 5.1]
Let $Y$ be a complete riemannian manifold
all ends of whose having infinite volume and such that 
$\lambda _{1}(Y) > 0$, then the natural morphism 
$$
H^{1}_{c}(Y,\Bbb R)
\to
H^{1}_{L^{2}}(Y,\Bbb R)
$$  
is injective.
In particular any $\alpha \in H^{1}_{c}(Y,\Bbb R)$ admits a representative
$\bar{\alpha}$ which is in $L^{2}(Y,\Bbb R)$.
\end{lemma}

\begin{cor}
Let $C'$ be as above and assume that there exists a compact 
essential hypersurface $Z' \subset \tilde{X}/C'$.
Then there exists an harmonic $n-1$-form $\omega$ in $L^{2}(\tilde{X}/C')$
such that $\int_{Z'}\omega \neq 0$. 
\end{cor}

{\bf Proof :}

\bigskip
\noindent
Let $\alpha \in H^{1}_{c}(\tilde{X}/C',\Bbb R)$ a Poincar\'{e} dual of 
$[Z'] \in H_{n-1}(\tilde{X}/C',\Bbb R)$.
By definition of $\alpha$, for any 
$\beta \in H^{n-1}(\tilde{X}/C', \Bbb R)$,
one has

\begin{equation}
\int_{Z'} \beta = \int_{\tilde{X}/C'}\beta \wedge \alpha,
\end{equation}

(\cite{bt} p.51, note that $\tilde{X}/C'$ has a " finite good cover").

After Lemma 5.4 , $\alpha$ admits a
non trivial harmonic representative $\bar{\alpha}$ in
$L^{2}(\tilde{X}/C')$. (In order to apply the Lemma 5.4, one has to check that all ends of 
$\tilde{X}/C'$ have infinite volume, ie
for a compact $K \subset \tilde{X}/C'$ each unbounded connected component of
$\tilde{X}/C' - K$ has infinite volume: this comes from the fact that the 
injectivity radius of  $\tilde{X}/C'$ is bounded below by the 
injectivity radius of $X = \tilde{X} /\Gamma$ and the sectional
curvature bounded above by $-1$.)
The $(n-1)$-harmonic form $\omega = *\bar{\alpha}$, where $*$ is the Hodge operator,
is in $L^{2}(\tilde{X}/C')$ and verifies after (5.17) 

\begin{equation}
\int_{Z'}\omega = \int_{\tilde{X}/C'}\omega \wedge \bar{\alpha} =
\int_{\tilde{X}/C'}\omega \wedge *\omega = ||\omega||^{2}_{L^{2}(\tilde{X}/C')}
\neq 0. 
\end{equation}

$\Box$

We can now prove the proposition 5.2.
Let us briefly describe the idea. 
We consider the iterates $F'^{k}$ of $F' : \tilde{X}/C' \to \tilde{X}/C'$.
As $F'$ is homotopic to the identity map, $F'^{k}(Z')$ is homologous to $Z'$
and if $\omega$ is the harmonic form of the corollary 5.5 we have 

\begin{equation}
\int_{Z'} (F'^{k})^{*}(\omega) = \int_{Z'}\omega  = a \neq 0.
\end{equation}

We don't know if $F'^{k}(Z')$ converges or stays in a compact subset of 
$\tilde{X}/C'$ but we will show that $F'^{k}(Z')$ cannot entirely diverge 
in $\tilde{X}/C'$ and that there exists a $z' \in Z'$ such that 
$F'^{k}(z')$ subconverges to a point $x \in \tilde{X}/C'$ with
$|Jac_{n-1}F'(x)| = 1$.

>From (5.19) one gets

\begin{equation}
0 < |a| =  |\int_{Z'} (F'^{k})^{*}(\omega)| 
\leq \int_{Z'} |Jac_{TZ'} F'^{k}(z)|. || \omega (F'^{k}(z)|| dz 
\end{equation}

where 
$$|Jac_{TZ'} F'^{k}(z)| =
 ||  DF'^{k}(z)(u_{1}) \wedge DF'^{k}(z)(u_{2})\wedge... \wedge DF'^{k}(z)(u_{n-1}) ||
$$
and $(u_{1},..,u_{n-1})$
is an orthonormal basis of $T_{z}(Z')$.

Let us define
$\mathcal{B} = \{ z\in Z',   |Jac F'^{k}_{TZ'}(z)|$  {\it does not converge to 0} $\}$.

For $z\in Z'$ we define the sequence $z_{k}$ by
$z_{0} = z$ and $z_{k} = F'(z_{k-1})=F'^{k}(z) \in \tilde{X}/C'$. 

\begin{lemma}
There exists $z\in \mathcal{B}$ and a subsequence $z_{k_{j}}$ such that 
 $z_{k_{j}}$ converges to a point $x  \in \tilde{X}/C'$ with 
$|Jac_{n-1}(x)| = 1.$ 
\end{lemma}

{\bf Proof :}

\bigskip
\noindent
We first remark that $\lim _{x\to \infty} ||\omega(x)|| = 0$. This follows 
the following facts: $\omega$ is harmonic, $\omega \in L^{2}(\tilde{X}/C')$ and 
the injectivity radius of $\tilde{X}/C'$ is bounded below by a positive constant.

Let us assume that for all $z\in \mathcal{B}$ the sequence $z_{k}$ diverges in
$\tilde{X}/C'$. Then we have, for all $z\in \mathcal{B}$,

\begin{equation}
||\omega (z_{k})|| = ||\omega (F'^{k} (z))|| \to 0
\end{equation}

whenever $k$ tends to $\infty$ because of the previous remark.

On the other hand, as $||\omega (F'^{k}(z))|| \leq  C$ and $|Jac _{n-1}F'^{k}| \leq 1$,
it follows from (5.20)

\bigskip
$$
\lim _{k\to \infty} |\int_{Z'} (F'^{k})^{*}(\omega)| \leq
$$
$$
\lim _{k\to \infty}[ \int_{\mathcal{B}} ||\omega(F'^{k}(z))||dz 
+  C \int_{Z'-\mathcal{B}} |Jac_{Z'} F'^{k}(z)|dz = 0
$$

which contradicts our assumption. 

Thus there exists a point $z\in Z'$ such that 
\begin{equation}
|Jac_{Z'} F'^{k}(z)| \to \alpha \neq 0
\end{equation}

and such that there exists a subsequence 
$z_{k_{j}}= F'^{k_{j}}(z)$ with

\begin{equation}
lim  _{j\to \infty}z_{k_{j}}= x \in \tilde{X}/C'.
\end{equation}

The property (5.22) comes from the fact that the sequence
$|Jac_{Z'} F'^{k}(z)|$ doesn't tend to zero and is decreasing (because 
$|Jac_{n-1} F'| \leq 1$).  

Let us define
$$
E_{0} = T_{z}Z' , E_{1}= DF'(z)(E_{0})
$$
 and 
$$
E_{k}= DF'(z_{k-1})(E_{k-1}) \subset T_{z_{k}}(\tilde{X}/C').
$$
As  $z_{k_{j}} \to x$ we can assume, after extracting again a subsequence, that
$E_{k_{j}} \to E \subset T_{x}(\tilde{X}/C')$.
On the other hand we also have
\begin{equation}
|Jac_{Z'}F'^{k}(z)|=
|Jac_{E_{k-1}}F'(z_{k-1})||Jac_{E_{k-2}}F'(z_{k-2})|...|Jac_{E_{0}}F'(z)|
\end{equation}

We know that 
$|Jac_{E_{k}} F'(z_{k})| = 1-\epsilon _{k}$ where $0 \leq\epsilon _{k} <1$.
As $z\in \mathcal{B}$, we have
$$
\lim _{k\to \infty} \pi _{j=1}^{k}(1-\epsilon _{j}) = \alpha > 0
$$
therefore $\lim _{k\to \infty}\epsilon _{k} = 0$ 
and by continuity we have
$|Jac_{E} F'(x)| = 1.$ $\Box$

Now we can finish the proof of the step 1.
Let consider a lift of $x$ in $\tilde{X}$ and 
$E$ in $T\tilde{X}$ that we again call $x$ and $E$.
Then we have 
$|Jac_{E}\tilde{F'}(x)| = 1$.
 $\Box$

Let us remark that corollary 4.4 and (5.19) give another proof of the inequality
$\delta \geq n-2$, which does not use the isosystolic inequality, ie. Theorem 3.3.

\bigskip
\noindent

{\bf Step 2 : The weak tangent of $\partial \tilde{X}$ and $\Lambda _{C}$}

\bigskip

We first recall the definition of the Gromov-distance on $\partial \tilde{X}$.

For two arbitrary points $\theta$ and $\theta'$ in 
$\partial \tilde{X}$ let us define 

\begin{equation}
l(\theta, \theta') = inf \{ t >0 / dist(\alpha_{\theta}(t) , \alpha _{\theta'}(t)) = 1\}
\end{equation}

and

\begin{equation}
d(\theta , \theta') = e^{-l(\theta, \theta') } 
\end{equation}

then $d$ is a distance on $\partial\tilde{X}$.

We now recall a few definitions following \cite{bk}.
A complete metric  space $(S,\bar{d})$ is a weak tangent of a metric space
$(Z,d)$ if there exist a point $0\in S$, a sequence of points 
$z_{k}\in Z$ and a sequence of positive
real numbers $\lambda_{k} \to\infty $
such that the sequence of pointed metric spaces 
$(Z,\lambda_{k}d,z_{k})$ converges in the pointed Gromov-Hausdorff topology to    
$(S,\bar{d},0)$ where $(Z,\lambda_{k}d)$ stands for the set $Z$ endowed with
the rescaled metric $\lambda_{k}d$.

Let us recall that the sequence of metric spaces
$(Z_{k},d_{k},z_{k})$  converges to
$(S,\bar{d},0)$ in the pointed Gromov-Hausdorff topology if the following 
conditions hold, (cf. (B-B-I), definition 8.1.1).

\begin{definition} We say that the sequence of metric
spaces $(Z_{k}, d_{k}, z_{k})$ converges to $(S , \bar{d}, o)$ if  
for any $R>0$, $\epsilon >0$ there exists $k_{0}$ such that
for any $k\geq k_{0}$ there exists a (non necessary continuous) map
 $f: B(z_{k},R) \rightarrow S$ such that 

\noindent
(i) $ f(z_{k}) = 0$,

\noindent
for any two points $x$ and $y$ in $B(z_{k},R)$,

\noindent
(ii)  $|\bar{d}(f(x),f(y)) - d_{k}(x,y)| \leq \epsilon$,

\noindent
and 

\noindent
(iii) the $\epsilon$-neighborhood of the set $f(B(z_{k},R)$ contains $B(0,R-\epsilon)$.
\end{definition}

In the previous definition, $B(z_{k},R)$ stands for the ball of radius $R$ centered at the point 
$z_{k}$ in $(Z_{k}, d_{k})$.

For a metric space $(Z,d)$ we will denote 
$WT(Z,d)$ the set of weak tangents of $(Z,d)$.

Let $\Gamma$ a cocompact group of isometries of $(\tilde{X},\tilde{g})$ a 
$n$-dimensional 
Cartan Hadamard manifold of sectional curvature $K_{\tilde{g}} \leq -1$
and $C'$ a subgroup of $\Gamma$.
The limit set $\Lambda_{\Gamma}$ of $\Gamma$ is the full boundary
$\partial \tilde{X}$ namely a topological $(n-1)-$dimensional sphere
$S^{n-1}$. We endow $\partial \tilde{X}$ with 
the Gromov distance $d$ defined in (5.26).
In \cite{bk} Lemma 5.2, M.Bonk and B.Kleiner show among other properties the following

\begin{lemma}
For
any weak tangent space $(S,\bar{d})$ in $WT((\partial \tilde{X},d) )$, 
$S$ is homeomorphic to $\partial \tilde{X}$ less a point, thus to 
$\Bbb{R}^{n-1}$.
\end{lemma}

In fact the crucial asumption in the above lemma, coming from the cocompacness of
$\Gamma$, is the property that any triple of points in $\partial \tilde{X}$ can be uniformly
separated by an element of $\Gamma$, ie. there is $\delta >0$ such that for any three points
$\theta_{1}, \theta_{2}, \theta_{3} \in \partial \tilde{X}$
there exists a $\gamma \in \Gamma$ such that $d(\gamma\theta_{i},\gamma \theta_{j}) \geq \delta$
for all $1\leq i \neq j \leq 3$.
Following the argument of M.Bonk and B.Kleiner one can show that if $C'$ is a subgroup of
$\Gamma$ such that one weak tangent $(S,\bar{d})$  of $(\Lambda_{C'},d)$ is in 
$WT(\partial \tilde{X},d)$ 
and enough triples of points of $\Lambda_{C'}$ can be uniformly separated by elements of $C'$,
then $S$ is homeomorphic to $\Lambda_{C'}$ less a point. In particular, 
$\Lambda_{C'}$ is homeomorphic to $\partial \tilde{X}$.

\begin{lemma}
Let $\mathcal{L} \subset \partial \tilde{X}$
be a closed $C'$-invariant set and 
$\theta_{0} \in \mathcal{L}$. We assume that there exist
a sequence of positive
real numbers $\lambda_{k} \to\infty $
such that the sequence of pointed metric spaces 
$(\mathcal{L},\lambda_{k}d,\theta_{0})$ converges in the pointed Gromov-Hausdorff topology to    
$(S,\bar{d},0)$ where $(S,\bar{d},0)$ is a weak tangent of
 $(\partial \tilde{X},d)$.
We also assume that there exist positive constants $C$ and $\delta$, a sequence of 
points $\theta_{0}^{k} = \theta_{0}, \theta_{1}^{k}, \theta_{2}^{k} \in \mathcal{L}$
and a sequence of elements $\gamma_{k} \in C'$ such that 
$C^{-1} \lambda_{k}^{-1}\leq d(\theta_{i}^{k}, \theta_{j}^{k}) \leq  C \lambda_{k}^{-1}$
and $d(\gamma_{k}\theta_{i}^{k},\gamma_{k} \theta_{j}^{k}) \geq \delta$
for all $0\leq i\neq j \leq 2$. Then, $S$ is homeomorphic to $\mathcal{L}$ less a point.
In particular $\mathcal{L}$ is homeomorphic to  $\partial \tilde{X}$.
 
\end{lemma}

The proof of this lemma is postponed in the Appendix.

\bigskip
\noindent
{\bf Step 3 : The limit set $\Lambda_{C'}$ of $C'$ 
and the limit set $\Lambda_{C}$ of $C$ are equal to a topological equator.}

\bigskip

We have shown in step 1 that $\Lambda_{C}$ is a subset of some 
topological equator $E(\infty)$.

Let $o\in \tilde{X}$ and $E \subset T_{o}\tilde{X}$ be such that
$|Jac_{E} \tilde{F}(o)| = 1$ and $E(\infty)$ is the equator associated to $E$. 

Recall that there exists a subgroup $C'$
of $C$ which globally preserves an hypersurface  
$\tilde{Z'} \subset \tilde{X}$ and that $\tilde{Z'}/C' \subset{\tilde{X}/C'}$ is compact.  
Furthermore $\tilde{Z'}$ separates $\tilde{X}$
into two connected components
$\tilde{U}$ et $\tilde{V}$.
We can assume that $\tilde{U}$ and $\tilde{V}$ are
 globally invariant by $C'$ after having replaced 
$C'$ by an index $2$ subgroup.

The limit set $\Lambda_{C'}$ of $C'$ is contained in $\Lambda_{C}$,
therefore $\Lambda_{C'}\subset E(\infty)$. 

We will show that $\Lambda_{C'} = E(\infty)$. 

For any subset $W \in \tilde{X}$ we define 
the boundary at infinity  $\partial W$ of $W$ by

\begin{equation}
\partial W = Cl(W) \cap \partial \tilde{X}
\end{equation}

where $Cl(W)$ stands for the closure of $W$ in $\tilde{X}\cup \partial \tilde{X}$.

As $\tilde{Z'}/C'$ is compact, $\tilde{Z'}$ is at bounded distance of the 
orbit $C'z$ of some point $z$ in $\tilde{Z'}$, thus

\begin{equation}
Cl(\tilde{Z'}) \cap \partial \tilde{X} = \Lambda_{C'}
\end{equation}

By definition we have 
$\Lambda _{C'} \subset \partial \tilde{U}$
and $\Lambda _{C'} \subset \partial \tilde{V}$.

\begin{lemma}
Let us assume that $\Lambda_{C'} \neq E(\infty)$, then either 
$\partial \tilde{U} = \Lambda _{C'}$ or 
$\partial \tilde{V}  = \Lambda _{C'}$.
\end{lemma}

{\bf Proof:}
We know that $\Lambda _{C'} \subset \partial \tilde{U}$
and $\Lambda _{C'} \subset \partial \tilde{V}$.

Let us assume that the conclusion of the lemma is not true, so
there is $\zeta \in \partial \tilde{U} - \Lambda_{C'}$ and
$\theta \in \partial \tilde{V} - \Lambda_{C'}$.

As $\Lambda _{C'} \subset E(\infty)$,
$\Lambda _{C'} \neq E(\infty)$ and $\partial \tilde{X}$ is a sphere,
any two points of $\partial \tilde{X} - \Lambda_{C'}$ can be joined by a 
continuous path contained in $\partial \tilde{X} - \Lambda_{C'}$ and so does
$\zeta$ and $\theta$, joined by such a path $\alpha$.

The set $Z' \cup \Lambda_{C'}$ is a closed subset of 
$\tilde{X} \cup \partial \tilde{X}$ thus there is an open connected neighborhood $W$
of $\alpha$ in  $\tilde{X} \cup \partial \tilde{X}$ contained in the complementary
of $Z' \cup \Lambda_{C'}$.

As $\zeta$ and $\theta$ can be approximated by points in $\tilde{U}$ and $\tilde{V}$
respectively there 
exist points $x\in \tilde{U} \cap W$ and $y\in \tilde{V} \cap W$ that can be joined by 
a continuous path by connectedness of $W$, which leads to a contradiction. $\Box$ 
\bigskip

\begin{rem}
In fact, we are going to show that under the assumption $\delta(C') = n-2$,
it is impossible to have $\partial \tilde{U}= \Lambda_{C'}$ or 
 $\partial \tilde{V}= \Lambda_{C'}$. 
\end{rem}

For any $x\in \tilde{X}$ and $\theta \in \partial \tilde{X}$ let us denote 
$HB(x,\theta)$ the open horoball centered at $\theta$ and passing through $x$.

\begin{lemma}
Let us assume that $\partial \tilde{U} = \Lambda_{C'}$. Then there exist 
$\theta_{0} \in \Lambda_{C'}$ and $z'\in \tilde{Z'}$
such that $HB(z',\theta_{0}) \subset \tilde{U}$.
\end{lemma}

{\bf Proof :}
Let us recall that $\tilde{X}/ C' - Z' = U \cup V$
where $U=\pi(\tilde{U})$ $V=\pi(\tilde{V})$ and $\pi : \tilde{X} \to \tilde{X} / C'$
is the projection .

We know that $U$ and $V$ are unbounded.
Let $x_{n}$ a sequence of points in $U$ such that 
$dist(x_{n}, Z') \to \infty$.
Let $z_{n}\in Z'$ such that $dist(x_{n}, Z') = dist(x_{n}, z_{n})$.
We consider a fundamental domain $D \subset \tilde{Z'}$ of $C'$.
There exist lifts $\tilde{z}_{n}\in D$ and $\tilde{x}_{n} \in \tilde{U}$ 
such that
$dist(\tilde{x}_{n},Z') = dist(\tilde{x}_{n},\tilde{z}_{n})$ tends
to infinity.

By compactness we can assume that 
a subsequence $\tilde{x}_{n_{j}}$ converges to a point 
$\theta_{0} \in \partial \tilde{X}$ and 
$\tilde{z}_{n_{j}}$ also converges to a point $\tilde{z} \in \bar{D}$.
Furthermore the sequence of open balls 
$B(\tilde{x}_{n_{j}}, dist( \tilde{x}_{n_{j}}, \tilde{z}_{n_{j}})) \subset \tilde{U}$
converges to the open horoball $HB(\theta_{0},\tilde{z}) \in \tilde{U}$. $\Box$

\begin{prop}

$\Lambda_{C'} = E(\infty)$ and $\Lambda_{C} = E(\infty)$.
\end{prop}

We first describe the idea of the proof and next state some facts we will need 
in order to do it.

As $\Lambda_{C'} \subset \Lambda_{C} \subset E(\infty)$ the proposition
boils down to proving that $\Lambda_{C'} = E(\infty)$.
Let us assume $\Lambda_{C'} \neq E(\infty)$ and find a contradiction.

We will show that for any sequence $\theta_{i}$ converging to 
$\theta_{0}$ the geodesic starting at 
a point $o$ such that $|Jac_{n-1}\tilde{F}(o)|=1$ and ending at $\theta_{i}$
crosses the hypersurface $\tilde{Z'}$ in  a point $z_{i}$.
For an {\bf appropriate choice} of such a sequence $\theta_{i}$
(roughly speaking, the sequence $\theta_{i}$ 
is chosen to be converging to $\theta_{0}$ ''transversally to $\Lambda_{C'}$''),
 the shadow 
(defined below) projected from
$o$ through some geodesic ball $B(z_{i}, r)$ will not intersect 
$\Lambda_{C'}$. 
On the other hand this shadow has to meet the limit set $\Lambda _{C'}$ 
because of the shadow lemma of D. Sullivan, which leads to a contradiction.

Precisely, by Lemma (5.10) and Lemma (5.12) we know that $\partial \tilde{U} = \Lambda _{C'}$
and that there exists an open horoball $HB(\theta_{0},\tilde{z})
 \subset \tilde{U}$
centered at a point $\theta \in \Lambda_{C'}$,
and whose closure contains a point $\tilde{z} \in \tilde{Z'}$.

Let $o\in \tilde{X}$ and $E\in T_{o}\tilde{X}$ an hyperplane such that
$\tilde{F}(o) = o$, $|Jac_{E}\tilde F(o)| = 1$, and 
$\Lambda_{C'} \subset E(\infty)$ where $E(\infty)$ is the topological 
equator associated to $E$.

For each $\theta \in \partial \tilde{X}$ we denote by 
$\alpha_{\theta}$ the geodesic starting from $o$ and such that 
$\alpha_{\theta}(+\infty) = \theta$.

Let $\theta_{i} \in \partial \tilde{X} - E(\infty) =\partial \tilde{V}-\partial \tilde{U}$
 be a
 sequence converging to $\theta_{0}$.
By continuity, for each $i$ large enough, the geodesic $\alpha_{\theta_{i}}$
spends some time 
inside the horoball  
$HB(\theta_{0},\tilde{z}) \subset \tilde{U}$ 
and ends up inside $\tilde{V}$ because $\theta_{i}$ converges to 
$\theta_{0}$ and $\theta_{i}$ belongs to
 $\partial \tilde{V} - \partial \tilde{U}$.

Thus $\alpha_{\theta_{i}}$ eventually crosses $\tilde{Z'}$.
Let $z_{i} \in \alpha_{\theta_{i}} \cap \tilde{Z'}$.
As $\tilde{Z'}/C'$ is compact, there is an element $\gamma_{i} \in C'$
such that $z_{i} = \gamma _{i}(x_{i})$ where $x_{i}$ is a point in the closure
$\bar{D}$ of a fundamental domain $D$ for the action of $C'$ on $\tilde{Z'}$.
The points $\gamma _{i}(x_{i})$ and  $\gamma _{i}(o)$ stay at bounded distance 
because
$dist(\gamma _{i}(x_{i}),\gamma _{i}(o)) = dist(x_{i}, o) \leq dist( o, D) + diam D$.
In particular, 
$lim_{i\to \infty} \gamma_{i}(o) = \theta_{0}$.

We have proved the 

\begin{lemma}
Let $\theta_{i} \in \partial \tilde{X}- E(\infty)$ be a sequence
 which converges to $\theta_{0}$.
There exists a constant $A$ such that
for $i$ large enough there exists $z_{i} \in \tilde{Z'}
 \cap \alpha_{\theta_{i}}$
and $\gamma_{i} \in C'$ such that $dist(z_{i},\gamma_{i}(o)) \leq A$ and both
$z_{i}$ and $\gamma_{i}(o)$ converge to $\theta_{0}$. 
\end{lemma}

Let $x$ and $y$ two points in $\tilde{X}$.

 We define the shadow $\mathcal{O}(x,y,R)\subset \partial \tilde X$ 
of the ball $B(y,R)$ enlighted from the point $x$ by

\begin{equation}
\mathcal{O}(x,y,R) = \{ \alpha(+\infty) \}
\end{equation}
 
where $\alpha$ runs through the set of geodesic rays starting from $x$ and meeting 
$B(y , R)$.

Let $\{\mu_{x}\}_{x}$ be a family of Patterson measures associated
 to the discrete 
group $C'$ with critical exponent $\delta' = \delta(C')$.

The following shadow lemma is due to D.Sullivan.

\begin{lemma} \cite{su1}, \cite{ro}, \cite{yue}.
There exist positive constants $C$ and $R$ such that
for any $y$ in $\tilde{X}$,
$\nu _{y}(\mathcal{O}(y, \gamma(y),R)) \geq C e^{\delta' d(y,\gamma(y))}$
\end{lemma}

\begin{cor}
Let $z_{i}$ be defined in lemma (5.14), then we have
$\mathcal{O}(o, z_{i}, R+A) \cap \Lambda_{C'} \neq \emptyset$ for $i$ large enough.
\end{cor}

We now prove that for a good choice of $\theta_{i}$, the shadow 
$\mathcal{O}(o, z_{i}, R+A) $ (with $z_{i}$ associated to $\theta_{i}$ as
in lemma 5.14)
never meet $\Lambda_{C'}$ for all large $i$'s, ie. for any $\theta \in \Lambda_{C'}$
the geodesic $\alpha_{\theta}$ does not cross 
$B(z_{i}, R+A)$.
We have no control on the radius $R$ coming from the shadow lemma nor on
the constant $A$ but we will show 

\begin{prop}
There exists a sequence 
$\theta_{i} \in \partial \tilde{X} - \Lambda_{C'}$
such that 
$\theta_{i}$ converges to $\theta_{0}$ 
and 
$$
lim_{i \to \infty} inf _{\theta\in \Lambda_{C'}} dist(z_{i}, \alpha _{\theta}) = + \infty
$$
where  
$z_{i} = \tilde{Z}' \cap \alpha_{\theta_{i}}$ has been constructed in
 lemma (5.14).
\end{prop}

\begin{cor}
For $i$ large enough,
$\mathcal{O}(o, z_{i}, R+A) \cap \Lambda_{C'} = \emptyset$.
\end{cor}

The corollary (5.16) and the corollary (5.18)
 lead to a contradiction, which ends the proof
of the proposition (5.13).

The end of the paragraph is devoted to proving the proposition (5.17).

\begin{lemma}
Let $\theta_{i}$ be a sequence of points in $\partial \tilde{X}$ converging 
to $\theta_{0}$ and $z_{i}$ constructed in lemma (5.14). Assume that 

\noindent
$liminf_{i \to \infty} inf _{\theta\in \Lambda_{C'}} dist(z_{i},
 \alpha _{\theta}) = C <+ \infty$
then 
$lim _{i\to \infty} \frac{d(\theta_{i},\Lambda_{C'})}
{d(\theta_{i},\theta_{0})} = 0$.
\end{lemma}

{\bf Proof :} 
We first show that 

\begin{equation}
\lim _{i\to \infty} dist(z_{i},\alpha_{\theta_{0}}) = \infty
\end{equation}

Recall that, for any $z \in \tilde{X}$ and $\theta \in \partial \tilde{X}$,
$B(z,\theta)$ equals the decreasing limit as $t$ tends to infinity 
of $dist(z, \alpha_{\theta}(t)) - dist(o, \alpha_{\theta}(t))$ 
where $\alpha_{\theta}(t)$ is the geodesic ray 
joigning $o$ to $\theta$. Therefore, as the points $z_{i}\in \tilde{Z}$ belongs 
to the complementary of the fixed horoball $HB(\tilde{z},\theta_{0})$,
we have,

\begin{equation}
dist(z_{i}, \alpha_{\theta_{0}}(T_{i}) \geq T_{i} + B(\tilde{z},\theta_{0})
\end{equation}
 
where $dist(z_{i}, \alpha_{\theta_{0}}(T_{i})) = dist(z_{i},\alpha_{\theta_{0}})$.

On the other hand, as
$z_{i}$ tends to $\theta_{0}$, $T_{i}$ tends to infinity so (5.30) is proven.

Let $t_{i}$ be such that $z_{i} = \alpha_{\theta_{i}}(t_{i})$. 
By (5.30) we have 

\begin{equation}
\lim _{i\to \infty} dist(\alpha_{\theta_{i}}(t_{i}) ,\alpha_{\theta_{0}}(t_{i})) =\infty.
\end{equation}

Let $u_{i}$ be such that 

\begin{equation}
\rm dist(\alpha_{\theta_{i}}(u_{i}), \alpha_{\theta_{0}}(u_{i})) = 1,
\end{equation}

then in particular $u_{i}\leq t_{i}$ for $i$ large enough and
by the triangle inequality we have

\begin{equation}
dist(\alpha_{\theta_{i}}(t_{i}), \alpha_{\theta_{0}}(t_{i})) \leq 2(t_{i} - u_{i}) +1.
\end{equation}

By (5.32) we get 

\begin{equation}
\lim_{i\to \infty} (t_{i} - u_{i}) = +\infty.
\end{equation}

Let us assume there exists a sequence $\theta_{i}' \in \Lambda_{C'}$ and 
a constant $C$ such that 
\begin{equation}
 dist(z_{i},\alpha_{\theta_{i}'}) \leq C < +\infty.
\end{equation}

We can assume that $C \geq 1$.

Let $v_{i}$ be such that 

\begin{equation}
dist(z_{i},\alpha_{\theta_{i}'}) = dist(z_{i},\alpha_{\theta_{i}'}(v_{i})).
\end{equation}

By triangle inequality,

\begin{equation}
|t_{i} - v_{i}| \leq C
\end{equation}

and,

\begin{equation}
dist(\alpha_{\theta'_{i}}(t_{i}), \alpha_{\theta_{i}}(t_{i})) \leq 2C.
\end{equation}

On the other hand, as the curvature of $\tilde{X}$ is bounded above by $-1$,
a classical comparison theorem gives for any $t\in [0,t_{i}]$,

\begin{equation}
sinh (\frac{dist(\alpha_{\theta'_{i}}(t), \alpha_{\theta_{i}}(t))}{2}) \leq
sinh C. \frac{sinht}{sinht_{i}}. 
\end{equation}

Let $s_{i}$ be such that 

\begin{equation}
dist(\alpha_{\theta'_{i}}(s_{i}), \alpha_{\theta_{i}}(s_{i})) = 1.
\end{equation}

There are two cases. Either $s_{i} \geq  t_{i}$ or $s_{i} < t_{i}$.
If $s_{i} < t_{i}$, we get from
(5.39) and (5.40) the existence of a constant $A$ such that for any $i$,

\begin{equation}
s_{i} \geq  t_{i} - A,
\end{equation}

and this inequality also holds when $s_{i} \geq t_{i}$.

>From (5.42) we get 

\begin{equation}
\frac{d(\theta_{i}, \theta'_{i})}{d(\theta_{i},  \theta_{0})}
= e^{-s_{i}+ u_{i}}
\leq e^{A}e^{-t_{i}+u_{i}},
\end{equation}

therefore, thanks to (5.35) we obtain

\begin{equation}
lim_{i\to \infty}\frac{d(\theta_{i}, \theta'_{i})}{d(\theta_{i},  \theta_{0})}= 0.
\end{equation}

which ends the proof of lemma (5.19). $\Box$

\begin{lemma}
Let us assume that for every sequence $\theta_{i}$ of
 points in $\partial \tilde{X}$ converging 
to $\theta_{0}$,
$lim_{i \to \infty}
 \frac{d(\theta_{i},\Lambda_{C'})}{d(\theta_{i},\theta_{0})}=0$.
Let $\lambda_{k}\to\infty$ be such that the sequence of spaces
$(\partial \tilde{X},\lambda_{k}d,\theta_{0})$
converges to a space $(S,\bar{d},0)$ in the pointed 
Gromov-Hausdorff topology,
then the sequence of spaces $(\Lambda_{C'},\lambda_{k}d,\theta_{0})$
also converges to $(S,\bar{d},0)$.
\end{lemma}

{\bf Proof :} Let us define
$$
r(\epsilon) =: sup\{ \frac{d(\theta,\Lambda_{C'})}{d(\theta,\theta_{0})}
,\theta \neq \theta_{0} , d(\theta,\theta_{0})\leq \epsilon \}.
$$

The assumption says
that
 
\begin{equation}
lim_{\epsilon \to 0} r(\epsilon) = 0.
\end{equation}

For an arbitrary metric space $(Y,d)$ and $Y'$ a 
subset of $Y$, let us denote
$B_{(Y,d)}(y,R)$ the closed ball of $(Y,d)$ 
of radius $R$ centered at $y\in Y$, and
$\mathcal{U}_{\epsilon}^{(Y,d)}(Y')$ the $\epsilon$-neighborhood of $Y'$ in $(Y,d)$.
For a metric space $(Y,d)$ and a positive number $\lambda$, let us denote
$\lambda Y$ the rescaled space $(Y,\lambda d)$.

By definition of the function $r$, we have for any $R$,

$$
B_{\lambda_{k}\partial \tilde{X}}(\theta_{0},R)
\subset 
\mathcal{U}_{\epsilon_{k}}^{\lambda_{k}\partial\tilde{X}}
B_{\lambda_{k}\Lambda_{C'}}(\theta_{0},R+\epsilon_{k})
$$
\begin{equation}
\subset
B_{\lambda_{k}\partial \tilde{X}}(\theta_{0},R+2\epsilon_{k}).
\end{equation}

where $\epsilon_{k} =: R r(R/\lambda_{k})$.

Let us fix $\alpha >0$.
By definition 5.7, for any $R>0$, $\epsilon >0$, there exist a map
$f: B_{\lambda_{k}\partial \tilde{X}}(\theta_{0},R+\alpha) \rightarrow S$ such that 
for $k\geq k_{0}$,

\noindent
(i) $ f(\theta_{0}) = 0$,

\noindent
for any two points $x$ and $y$ in $B_{\lambda_{k}\partial\tilde{X}}(\theta_{0},R+\alpha)$ ,

\noindent
(ii)  $|\bar{d}(f(x),f(y)) - \lambda_{k}d(x,y)| \leq \epsilon$,

\noindent
and 

\noindent
(iii) $B_{(S,\bar{d})}(0,R +\alpha-\epsilon)\subset 
\mathcal{U}_{\epsilon}^{(S,\bar{d})}f(B_{\lambda_{k}\partial\tilde{X}}(\theta_{0},R+\alpha))$.

\bigskip
Moreover let us prove:

\noindent
(iv) $B_{(S,\bar{d})}(0,R -2\epsilon)\subset 
\mathcal{U}_{\epsilon}^{(S,\delta)}f(B_{\lambda_{k}\partial\tilde{X}}(\theta_{0},R))$.

Indeed, let $z\in B_{(S,\bar{d})}(0,R-2\epsilon)$. By (iii), there exists 
$\theta \in B_{\lambda_{k}\partial \tilde{X}}(\theta_{0}, R+\alpha)$ such that
$\bar{d}(z,f(\theta)) \leq \epsilon$. Since $\bar{d}(z,0) \leq R-2\epsilon$,
we thus deduce from triangle inequality $\bar{d}(f(\theta),0) \leq R-\epsilon$,
and therefore we get, thanks to (i) and (ii), $\lambda_{k}d(\theta, \theta_{0}) \leq R$. $\Box$

By (5.45),
for $\epsilon$ small enough, there exists $k_{1} \geq k_{0}$ such that for any $k \geq k_{1}$,
then $2\epsilon_{k} \leq \epsilon$ and 

$$
B_{\lambda_{k}\partial\tilde{X}}(\theta_{0},R+2\epsilon_{k}) \subset 
B_{\lambda_{k}\partial\tilde{X}}(\theta_{0},R+\alpha).  
$$

Therefore, by (5.46) and the above properties (i),(ii),(iii) and (iv) of the map $f$,
and the triangle inequality we get,

$$
B_{(S,\bar{d})}(0,R-2\epsilon) \subset \mathcal{U}_{\epsilon}^{(S,\bar{d})}
f (B_{\lambda_{k}\partial\tilde{X}}(\theta_{0},R)) 
$$

$$
\subset \mathcal{U}_{\epsilon}^{(S,\bar{d})}
f(\mathcal{U}_{\epsilon_{k}}^{\lambda_{k}\partial\tilde{X}}
B_{\lambda_{k}\Lambda_{C'}}(\theta_{0},R+\epsilon_{k}))
$$

$$
\subset \mathcal{U}_{2\epsilon + \epsilon_{k}}^{(S,\bar{d})}
f (B_{\lambda_{k} \Lambda_{C'}}(\theta_{0},R+\epsilon_{k})). 
$$

About the second inclusion above let us remark that the set
$\mathcal{U}_{\epsilon_{k}}^{\lambda_{k}\partial \tilde{X}}
B_{\lambda_{k}\Lambda_{C'}}(\theta_{0},R+\epsilon_{k})$
is contained in 
$B_{\lambda_{k}\partial \tilde{X}}(\theta_{0}, R+2\epsilon_{k})
\subset B_{\lambda \partial \tilde{X}}(\theta_{0}, R + \alpha)$,
so that we can apply $f$ to this set.

>From the above inclusions we obtain

$$
B_{(S,\bar{d})}(0,R-3\epsilon)
\subset \mathcal{U}_{2\epsilon +\epsilon_{k}}^{(S,\bar{d})}
f (B_{\lambda_{k} \Lambda_{C'}}(\theta_{0}, R))
$$
 
which implies the convergence of 
$(\Lambda_{C'}, \lambda_{k}d, \theta_{0})$ to $(S,\bar{d},0)$.

 $\Box$

\begin{cor}
Let us assume that for every sequence $\theta_{i}$ of
 points in $\partial \tilde{X}$ converging 
to $\theta_{0}$ and $z_{i}$ the sequence of points constructed 
in lemma (5.14), 
$liminf_{i \to \infty} inf _{\theta\in \Lambda_{C'}}
 dist(z_{i}, \alpha _{\theta})  <+ \infty$.
Let $\lambda_{k}\to \infty$ be such that the sequence of spaces
$(\partial \tilde{X},\lambda_{k}d,\theta_{0})$
converges to the space $(S,\delta,0)$ in the pointed 
Gromov-Hausdorff topology,
then the sequence of spaces $(\Lambda_{C'},\lambda_{k}d,\theta_{0})$
also converges to $(S,\delta,0)$.
\end{cor}

We will show now that there exist a sequence of 
points $\theta_{1}^{k}, \theta_{2}^{k} \in \Lambda_{C'}$ converging to 
$\theta_{0}$, such that the mutual distances 
$d(\theta_{1}^{k}, \theta_{2}^{k})$,
$d(\theta_{1}^{k}, \theta_{0})$,
$d(\theta_{2}^{k}, \theta_{0})$ is tending to zero at the same rate, and the 
triple $\theta_{1}^{k}, \theta_{2}^{k}, \theta_{0}$ can be uniformly separated
by elements $\gamma_{k} \in C'$. 

\begin{lemma}
Assume that every weak tangent of $(\partial \tilde{X},d)$ at $\theta_{0}$
 belongs to $WT(\Lambda_{C'},d)$, then
there exist positive constants $c$, $\delta$, a sequence $\epsilon_{k}$ tending to $0$
when $k$ tends to $\infty$, a sequence $\gamma_{k} \in C'$,
a sequence of points $\theta_{1}^{k}, \theta_{2}^{k} \in \Lambda_{C'}$
 such that for $i = 1, 2$,
 
\noindent
$c^{-1}\epsilon_{k}\leq d(\theta_{1}^{k},
 \theta_{2}^{k})\leq c \epsilon_{k}$,

\noindent
$c^{-1}\epsilon_{k}\leq d(\theta_{i}^{k}, \theta_{0}))\leq c \epsilon_{k}$
and

\noindent
$ d(\gamma_{k}\theta_{1}^{k}, \gamma_{k}\theta_{2}^{k})\geq \delta $,
$ d(\gamma_{k}\theta_{i}^{k}, \gamma_{k}\theta_{0}))\geq \delta$.

\end{lemma}

{\bf Proof :}
For any $x\in \tilde{X} \cup \partial \tilde{X}$
 and $y\in \tilde{X} \cup \partial \tilde{X}$ let us define
$\alpha_{x,y}$ the geodesic ray joining $x$ and $y$.
 Let $o\in \tilde{X}$ and $E \in T_{o}\tilde{X}$ be such that 
$|Jac_{E}\tilde{F}(o)| = 1$ and $E(\infty)$ the equator associated to $E$.
Let $\gamma_{k} \in C'$ be a sequence such that $\gamma_{k}(o)$ converges to 
the point $\theta_{0}$ where $\theta_{0} \in \Lambda_{C'}$ is the point coming from lemma 5.12.
In particular, according to that lemma, there exist a point $\tilde{z} \in \tilde{Z'}$ 
such that the hypersurface $\tilde{Z'}$ is contained in the complementary of the open horoball
$HB(\tilde{z}, \theta_{0})$.
We define  
$D : = dist(\tilde{z},o)$.
As $\tilde{Z'}$ lies outside the open horoball $HB(\tilde{z},\theta_{0})$, the points
$\gamma_{k}(o)$ belong to the complementary of the open horoball
$HB(\alpha_{\tilde{z},\theta_{0}}(D) ,\theta_{0})$.
By standard triangle comparison argument (comparison with the hyperbolic case) the angle  
$Angle( \alpha_{\gamma_{k}(o) , \theta_{0}} , \alpha_{\gamma_{k}(o) , o}) $
between the two geodesic rays
$\alpha_{\gamma_{k}(o) , \theta_{0}}$ and $ \alpha_{\gamma_{k}(o) , o}$ 
satisfies :
\begin{equation}
 lim _{k\to\infty} Angle( \alpha_{\gamma_{k}(o) , \theta_{0}} , \alpha_{\gamma_{k}(o) , o}) = 0.
\end{equation}

By equivariance we have $\Lambda_{C'} \subset (\gamma_{k}E)(\infty)$ where
$\gamma_{k}E \subset T_{\gamma_{k}(o)}\tilde{X}$.
For any $v\in T\tilde{X}$ let $\alpha_{v}$ be the geodesic ray such that 
$\dot{\alpha}_{v}(0) = v$.
Let us denote by $u_{k}$ the unit vector in $\gamma_{k}E$ such that 
$\alpha_{u_{k}}(+\infty) = \theta_{0}$ and let us choose some $w_{k}\in \gamma_{k}E$
such that $<u_{k},w_{k}> = 0$ (this is possible because $n-1\geq 2$.

We claim now that there exist $v_{k}\in \gamma_{k}E$ such that the angle
between $v_{k}$ and $w_{k}$ is not too far from $0$ or $\pi$, namely
\begin{equation}
|<v_{k},w_{k}>| \geq \frac{1}{(n-1)^{1/2}},
\end{equation}
and $\alpha_{v_{k}}(+\infty) \in \Lambda_{C'}$ or 
 $\alpha_{v_{k}}(-\infty) \in \Lambda_{C'}$.

Let us prove this claim.

According to Proposition 5.1 and to (5.11), the restriction to 
$\gamma_{k}E$ of the quadratic form
$h(u) = \int DB(\gamma_{k}(o), \theta)(u)^{2} d\mu_{\gamma_{k}(o)}(\theta)$
 verifies
\begin{equation}
h_{\gamma_{k}E}(u) = \frac{||u||^{2}}{n-1}.
\end{equation}

Therefore, if for all $u\in \gamma_{k}(E)$  
such that 
$\alpha_{u}(+\infty) = \theta \in \Lambda_{C'}$ we had $|<u, w_{k}>| < \frac{1}{(n-1)^{1/2}}$,
then one would get  
$h(w_{k}) < \frac{1}{n-1}$, which contradicts (5.49) and proves the claim.

In particular the angle between 
$u_{k}$ and $v_{k}$ is not too far from $\pi / 2$ for
$k$ large enough, ie.

\begin{equation}
|<u_{k}, v_{k}>| \leq \Big( \frac{n-2}{n-1}\Big)^{1/2},
\end{equation} 

and thanks to (5.47), we have for $k$ large enough

\begin{equation}
|<\dot{\alpha}_{\gamma_{k}(0),o}(0), v_{k}>|
 \leq \Big(\frac{n-\frac{3}{2}}{n-1}\Big)^{1/2}.
\end{equation}

Let us now assume for example that
$\theta_{k} = \alpha_{v_{k}}(+\infty) \in \Lambda_{C'}$.
Let us show that

\begin{equation}
lim _{k\to \infty }d(\theta_{0},\theta_{k}) = 0.
\end{equation}

Assume that (5.52) is not true. Then, one can assume 
after extracting a subsequence that $\theta_{k}$ converges to $\theta \neq \theta_{0}$. 
 Therefore the geodesic rays 
$\alpha_{\gamma_{k}(0),o}$ and $\alpha_{v_{k}}$ would converge to 
the geodesics $\alpha_{\theta_{0},o}$ and 
$\alpha_{\theta_{0},\theta}$ and thus the angle 
$Angle(\alpha_{\gamma_{k}(0),o},\alpha_{v_{k}})$
would converge to $0$. But this would contradict (5.51).

Let us now denote $\epsilon_{k} =: d(\theta_{k}, \theta_{0})$. According to 
(5.52), $lim_{k\to \infty}\epsilon_{k} = 0$.
We now consider the following sequence of pointed metric space 
$(\partial \tilde{X}, \epsilon_{k}^{-1}d, \theta_{0})$, a subsequence of which being
converging to some metric space $(S,\delta)$, cf[].
 For convenience we still denote by the same index 
$k$ the subsequence. By the corollary 5.21, the sequence 
 $(\Lambda_{C'}, \epsilon_{k}^{-1}d, \theta_{0})$  also converges to $(S,\delta)$. 
According to lemma 5.8, the space $S$ is homeomorphic to $\Bbb{R}^{n-1}$. In particular there exist
a sequence of points $\theta _{k}' \in \Lambda_{C'}$ and a constant $c$ such that 

\begin{equation}
c^{-1}\epsilon_{k}\leq d(\theta_{k},\theta_{k}')\leq c \epsilon_{k},
\end{equation}

\begin{equation}
c^{-1}\epsilon_{k}\leq d(\theta_{k}', \theta_{0}))\leq c \epsilon_{k}.
\end{equation}

Thus, the points $\theta_{1}^{k}=\theta_{k}$ and $\theta_{2}^{k} = \theta_{k}'$
satisfy the two first properties of lemma 5.22.

In order to complete the proof of lemma 5.22, we will show that 
the elements $\eta _{k} =:\gamma_{k}^{-1}$ uniformly separate $\theta_{0}$,
$\theta_{1}^{k}$ and $\theta_{2}^{k}$.

Thanks to (5.50) the angle at $\gamma_{k}(o)$ between $\theta_{1}^{k}$ and 
$\theta_{0}$ is uniformly bounded away from $0$ and $\pi$ and so does 
the angle at $o$ between $\gamma_{k}^{-1}(\theta_{1}^{k})$
and $\gamma_{k}^{-1}(\theta_{0})$.
Therefore, as the angle is H\"{o}lder-equivalent to the distance $d$, 
 cf. \cite{ka}, there is a constant $c$
such that 
\begin{equation}
d(\gamma_{k}^{-1}(\theta_{1}^{k}),\gamma_{k}^{-1}(\theta_{0})) \geq c.
\end{equation}

Now the cocompact group $\Gamma$ acts uniformly quasi-conformally on 
$(\partial \tilde{X},d)$, (\cite{bk} and \cite{va} Theorem 5.2), and so does $C' \subset  \Gamma$, 
therefore 

\begin{equation}
d(\gamma_{k}^{-1}(\theta_{1}^{k}),\gamma_{k}^{-1}(\theta_{2}^{k})) \geq c,
\end{equation}

and

\begin{equation}
d(\gamma_{k}^{-1}(\theta_{2}^{k}),\gamma_{k}^{-1}(\theta_{0})) \geq c.
\end{equation}

which ends the proof of lemma 5.22. $\Box$

{\bf Proof of Proposition 5.17 :}

Let us assume that for every sequence $\theta_{i}$ of
 points in $\partial \tilde{X}$ converging 
to $\theta_{0}$, 
$liminf_{i \to \infty} inf _{\theta\in \Lambda_{C'}}
 dist(z_{i}, \alpha _{\theta})  <+ \infty$,
then by corollary 5.21 and lemma 5.22 there exist 
 a positive constant $c$, a sequence $\epsilon_{k}$ tending to $0$
when $k$ tends to $\infty$, a sequence $\gamma_{k} \in C'$,
a sequence of points $\theta_{1}^{k}, \theta_{2}^{k} \in \Lambda_{C'}$
 such that
 
\noindent
$c^{-1}\epsilon_{k}\leq d(\theta_{1}^{k},
 \theta_{2}^{k})\leq c \epsilon_{k}$,

\noindent
$c^{-1}\epsilon_{k}\leq d(\theta_{i}^{k}, \theta_{0}))\leq c \epsilon_{k}$
and 

\noindent
$ d(\gamma_{k}\theta_{1}^{k}, \gamma_{k}\theta_{2}^{k})\geq \delta $,
$ d(\gamma_{k}\theta_{i}^{k}, \gamma_{k}\theta_{0}))\geq \delta$.

Applying lemma 5.9 for $\mathcal{L}= \Lambda_{C'}$ and $\lambda_{k} = \epsilon _{k}^{-1}$ we conclude
that $\Lambda_{C'}$ is homeomorphic to $\partial \tilde{X}$, which is 
impossible because $\Lambda _{C'}$ is contained in a topological equator 
$E(\infty)$. $\Box$

\bigskip
{\bf Step 4 : $C'$ and $C$ are convex cocompact.}

\bigskip

We first define convex cocompactness. For a discrete group $C$ of isometries 
acting on a Cartan Hadamard manifold of negative sectional curvature with 
limit set $\Lambda _{C}$, one defines the geodesic hull
 $\mathcal{G}(\Lambda_{C})$ of
$\Lambda _{C}$ as the set of all geodesics both ends of whose belong to 
 $\Lambda _{C}$.

The geodesic hull of  $\Lambda _{C}$ is a $C$ invariant set. 
One says that $C$ is convex cocompact if  $\mathcal{G}(\Lambda_{C})/C$
is compact.

\begin{lemma}
C' is convex cocompact.
\end{lemma}

{\bf Proof :} Let us denote $\pi : \tilde{X} \to \tilde{X}/C'$ the projection.
Assume that $C'$ is not convex cocompact.
Then, there exist a sequence $x_{n} \in \mathcal{G}(C')$ such that $x_{n}$
tends to infinity. In particular $dist(x_{n},Z') \to +\infty$,
where $Z'= \tilde{Z'}/C'$
 is the compact hypersurface which separates $\tilde{X}/C'$
in two unbounded connected components.
There exist lifts $\tilde{x}_{n}$ of $x_{n}$ such that 

\begin{equation}
\tilde{x}_{n} \to \theta_{0} \in \Lambda_{C'} 
\end{equation}

\begin{equation}
dist(\tilde{x}_{n}, \tilde{Z'})= dist(\tilde{x}_{n}, \tilde{z}_{n})
\end{equation}

where $\tilde{z}_{n} \in \tilde{Z'}$ is bounded.
Therefore there exist $\tilde{z} \in \tilde{Z'}$

such that 
$HB(\tilde{z},\theta_{0}) \subset \tilde{U}$, where $\tilde{U}$ is one
of the two connected components of $\tilde{X} -\tilde{Z'}$, the other being $\tilde{V}$.

We recall that $\mathcal{M}$ , $\mathcal{N}$ are the two connected components
of $\partial \tilde{X} - \Lambda_{C'}$. 
We also have $\partial \tilde{Z'} = \Lambda_{C'} = E(\infty)$, and after possibly replacing $C'$
by an index two subgroup, we can assume that $C'$ preserves $\tilde{U}$ and $\tilde{V}$. 

{\bf Claim :} There are the two following cases.

 Either one of the two boundaries 
 $\partial \tilde{U}$ or
$\partial \tilde{V}$
is equal to $\Lambda_{C'}$ (in this case the other boundary is 
equal to $\partial \tilde{X}$),
 or 
 $\partial \tilde{U} = \bar{\mathcal{M}}$
and $\partial \tilde{V} = \bar{\mathcal{N}}$, where $\bar{\mathcal{M}}$ and 
$\bar{\mathcal{N}}$ are the closure of $\mathcal{M}$ and $\mathcal{N}$.

Let us prove the claim.
We first remark that if there exist $\theta \in \partial \tilde{U} \cap \mathcal{M}$,
then $\mathcal{M} \subset \partial \tilde{U}$. 
Namely, let $\xi$ be any other point in $\mathcal{M}$ and 
$\alpha$ a continuous path in $\mathcal{M}$ joining
$\theta$ and $\xi$. Since the set $\tilde{Z'} \cup \Lambda_{C'}$
is a closed subset in $\tilde{X} \cup \partial \tilde{X}$, there exist an open connected
neighborhood $W$ of $\alpha$ in $\tilde{X} \cup \partial \tilde{X}$ contained in the complementary
of $\tilde{Z'} \cup \Lambda_{C'}$. Therefore, as $W \cap \tilde{U} \neq \emptyset$,
we have $ W \cap \tilde{X} \subset \tilde{U}$ and
 $\xi  \in \partial \tilde{U}$.
Let us assume that neither $\partial \tilde{U}$ nor $\partial \tilde{V}$ is
equal to $\Lambda_{C'}$. Then, each boundary
$\partial \tilde{U}$ and $\partial \tilde{V}$
contains  $\mathcal{M}$ or  $\mathcal{N}$.
But on the other hand,
since the set  $\tilde{Z'} \cup \Lambda_{C'}$ is closed,
 $(\partial \tilde{U} - \Lambda_{C'}) \cap (\partial \tilde{V}-\Lambda_{C'}) = \emptyset$
 thus we have 
 $\partial \tilde{U}=\mathcal{M}$ and $\partial \tilde{V}=\mathcal{N}$ or
the other way around and the claim is proved.

{\bf Case 1} : $\partial \tilde{U} = \Lambda_{C'}$ and 
$\partial \tilde{V} = \partial \tilde{X}$ or the other way around.

In this case, 
we are in the situation of the step 3, which leads to a contradiction, cf. remark 5.11.

{\bf Case 2} : $\partial \tilde{U}= \bar{\mathcal{M}}$ and $\partial \tilde{V} = \bar{\mathcal{N}}$.
 
In that case, assuming $C'$ is not convex-cocompact, there exist an open horoball
$HB(\theta_{0},\tilde{z}) \subset \tilde{U}$ where
 $\theta_{0} \in \Lambda_{C'}$, $\tilde{z} \in \tilde{Z'}$,
 $\partial \tilde{U} = \bar{\mathcal{M}}$ and $\partial \tilde{V}= \bar{\mathcal{N}}$.
We will find a contradiction in a similar way as in case 1, ie. step 3.
We consider a point $o \in \tilde{X}$ and an hyperplane $E \subset T_{o}\tilde{X}$ such that 
$|Jac_{E}\tilde{F}(o)| = 1$ and $\Lambda_{C'} = E(\infty)$.

Let $\theta_{i} \in \mathcal{N}$ be a sequence which converge to $\theta_{0}$.
By continuity, for $i$ large enough, the geodesic ray $\alpha_{o,\theta_{i}}$ spends some time
in $HB(\theta_{0},\tilde{z})\subset \tilde{U}$ and ends up in $\tilde{V}$ because $\theta_{i}$
converges to $\theta_{0}$ and $\theta_{i}$ belongs to 
$\mathcal{N} = \partial \tilde{V} - \Lambda_{C'}$. Therefore, $\alpha_{o,\theta_{i}}$ eventually
crosses $\tilde{Z'}$. Let $z_{i}$ be some point in $\tilde{Z'} \cap \alpha_{o,\theta_{i}}$.

We will prove the following Proposition, similar to the Proposition 5.18,

\begin{prop}
There exist a sequence $\theta_{i} \in \mathcal{N}$ such that $\theta_{i}$ converges
to $\theta_{0}$ and 
$$
lim_{i\to \infty} inf _{\theta \in \Lambda_{C'}} dist(z_{i}, \alpha_{\theta}) = + \infty
$$
where $z_{i} \in \tilde{Z'} \cap \alpha_{\theta_{i}}$.
\end{prop}

\begin{rem}
The difference between the propositions 5.24 and 5.18 is that we are looking for a sequence
$\theta_{i} \in \mathcal{N}$ instead of $\theta_{i} \in \partial \tilde{X} -\Lambda_{C'}$.
\end{rem}

Assuming the Proposition 5.24 we find a contradiction in the same way as in step 3.
Namely, as $\tilde{Z'}/C'$ is compact, the points $z_{i}\in \tilde{Z'} \cap \alpha_{\theta_{i}}$
stay at bounded distance from the $C'$-orbit of a fixed point, say, $o$, 
thus there exist a constant $A>0$ 
and elements $\gamma_{i} \in C'$ such that for any $i$, 

\begin{equation}
dist(z_{i},\gamma_{i} o)\leq A.
\end{equation}

>From (5.60) and the shadow lemma 5.15, we obtain
$\mathcal{O}(o,z_{i}, R+A) \cap \Lambda_{C'} \neq \emptyset$,
and on the other hand, from the proposition 5.24, we have
$\mathcal{O}(o,z_{i},R+A) \cap \Lambda_{C'} = \emptyset$,
which gives the contradiction.
It remains to prove the Proposition 5.24.

{\bf Proof of the proposition 5.24 : }
We argue by contradiction, like in the proof of the proposition 5.17.
Let us assume that there exist a constant $C>0$ such that  
for any sequence of points $\theta_{i} \in \mathcal{N}$ converging to $\theta_{0}$, 
$lim_{i\to \infty} inf _{\theta \in \Lambda_{C'}} dist(z_{i},\alpha_{\theta}) \leq C$,
then by lemma 5.19, we have for any such sequence $\theta_{i} \in \mathcal{N}$

\begin{equation}
lim_{i\to\infty} \frac{d(\theta_{i},\Lambda_{C'})}{d(\theta_{i},\theta_{0})} =0.
\end{equation}

The proof of the following lemma is the same as the proof of lemma 5.20.

\begin{lemma}
Let us assume that for any sequence $\theta_{i} \in \mathcal{N}$ conveging to $\theta_{0}$,
$lim_{i\to \infty} \frac{d(\theta_{i},\Lambda_{C'})}{d(\theta_{i},\theta_{0})} = 0$.
Let $\{\lambda_{k} \}$ be a sequence of positive numbers tending to $+\infty$
such that the sequence of spaces $(\partial \tilde{X}, \lambda_{k}d,\theta_{0})$
converges to a space $(S,\delta, 0)$ in the pointed Gromov-Hausdorff topology, then,
$(\mathcal{M},\lambda_{k}d,\theta_{0})$ also converges to $(S,\delta,0)$.
\end{lemma}

{\bf Proof : } Since $\Lambda_{C'} \subset \bar{\mathcal{M}}$, the assumption implies that
$lim _{\epsilon \to 0} r(\epsilon) =0$ where 
$$
r(\epsilon) = sup \{ \frac{d(\theta,\mathcal{M})}{d(\theta,\theta_{0})}, \theta \neq \theta_{0},
\theta \in \mathcal{N}, d(\theta,\theta_{0}) \leq \epsilon \}
$$

and the proof goes the same way as in lemma 5.20 replacing $\Lambda_{C'}$ by $\mathcal{M}$.
$\Box$

Similarly to the lemma 5.22, we have the 

\begin{lemma}
Let us assume that every weak tangent of $(\partial \tilde{X},d)$ at $\theta_{0}$ belongs
to $WT((\mathcal{M},d))$. There exist positive constant $c$, $\delta$, a sequence $\epsilon_{k}$ tending
to $0$ when $k$ tends to $+\infty$, a sequence of $\gamma_{k} \in C'$, a sequence of points
$\theta_{0}^{k}= \theta_{0},\theta_{1}^{k}, \theta_{2}^{k} \in  \mathcal{M}$
such that for $i \neq j \in \{0, 1, 2\}$, 

$c^{-1} \epsilon_{k} \leq d(\theta_{i}^{k}, \theta_{j}^{k}) \leq c \epsilon_{k}$ and 

$d(\gamma_{k}\theta_{i}^{k},\gamma_{k}\theta_{i}^{k}) \geq \delta.$
\end{lemma}

We can now end the proof of the proposition 5.24.
Let us assume that there exist a constant $C>0$ such that for every sequence $\theta_{i}$ of points in $\mathcal{N}$ 
converging to $\theta_{0}$,
$lim_{i\to \infty} inf_{\theta \in \Lambda_{C'}} dist(\theta_{i},\alpha_{\theta})\leq C$, then 
by (5.61), lemma 5.26 and lemma 5.27, there exist a sequence $\epsilon_{k}$ tending to $0$ when 
$k$ tends to $\infty$, a sequence $\gamma_{k} \in C'$, a sequence of points 
$\theta_{0}^{k}= \theta_{0},\theta_{1}^{k}, \theta_{2}^{k} \in  \mathcal{M}$
such that for $i \neq j \in \{0, 1, 2\}$, 

$c^{-1} \epsilon_{k} \leq d(\theta_{i}^{k}, \theta_{j}^{k}) \leq c \epsilon_{k}$ and 

$d(\gamma_{k}\theta_{i}^{k},\gamma_{k}\theta_{i}^{k}) \geq \delta.$

Applying the lemma 5.9 for $\mathcal{L} =\mathcal{M}$ and $\lambda_{k}= \epsilon_{k}^{-1}$, we 
conclude that $\mathcal{M}$ is homeomorphic to $\partial \tilde{X}$, which is impossible
because $\partial \tilde{X}$ is a sphere, and $\mathcal{M}$ is homeomorphic to 
an hemisphere.
This ends the proof of the proposition 5.24.
$\Box$

\begin{cor}
C is convex cocompact.
\end{cor}

{\bf Proof :} The subgroup $C'$ of $C$ is convex cocompact and the limit sets of 
$C'$ and $C$ coincide by step 3, therefore $C$ is convex cocompact. $\Box$.

\bigskip
{\bf Step 5:  C preserves a copy of the $(n-1)$-dimensional
hyperbolic space $\mathbb{H}^{n-1}$ totally geodesically embedded
in $\tilde{X}$.}

\bigskip
From the steps 1-4, we know that the groups $C$ and $C'$ are convex cocompact,
and that their limit set $\Lambda_{C}$ and $\Lambda_{C'}$ are equal to a topological
equator $E(\infty)$.

Let us consider the essential hypersurface $Z' \subset \tilde{X}/C'$.
We will show that there exist a minimizing current representing the class
of $Z'$ in $H_{n-1}(\tilde{X}/C',\mathbb{R})$ and that this minimizing current
lifts to a totally geodesic hypersurface embedded in $\tilde{X}$. We will then show that this
totally geodesic hypersurface is eventually hyperbolic.

We work in $\tilde{X}/C'$ and consider the essential hypersurface $Z'\subset \tilde{X}/C'$.
We will now prove that there exist a minimal current representing the class
of $Z'$ in $H_{n-1}(\tilde{X}/C',\mathbb{R})$.
Let $\{Z_{k}\}$ be a minimizing sequence of currents homologous to $Z'$.
The othogonal projection onto the convex core of $\tilde{X}/C'$ is distance 
nonincreasing and thus volume nonincreasing. Therefore we can assume that the $Z_{k}$'s are in the 
the convex core of $\tilde{X}/C'$, which is compact.  
By \cite{mo} (5.5), the sequence
 $\{Z_{k}\}$ subconverges to a minimal current $Z_{\infty}$
in $\tilde{X}/C'$.
By \cite{mo} (8.2), $Z_{\infty}$ is a manifold
with possible singularities of codimension greater than or equal to 8.
By corollary (4.4) and minimality we get that
 $|Jac_{n-1}\tilde{F}(x)| = 1$ at every regular points $x \in Z_{\infty}$.
We will use the fact that $|Jac_{n-1}\tilde{F}(x)| = 1$ at every regular points $x \in Z_{\infty}$
 in order to prove that $Z_{\infty}$ is 
a totally geodesic hypersurface.

\begin{lemma}
Let $x$ and $y$ two distinct points in $\tilde{X}$ and $E_{x}\subset T_{x}\tilde{X}$,
$E_{y}\subset T_{y}\tilde{X}$  be such that 
$Jac_{n-1}\tilde{F'}(x) = Jac_{E_{x}}\tilde{F'}(x) =1$ and
$Jac_{n-1}\tilde{F'}(y) = Jac_{E_{y}}\tilde{F'}(y) =1$. Then, 
the geodesic $\alpha_{x,y}$ (resp. $\alpha_{y,x}$) joining
$x$ and $y$ (resp. $y$ and $x$) satisfies
$\dot{\alpha}_{x,y}(0)\in E_{x}$, (resp. $\dot{\alpha}_{y,x}(0)\in E_{y}$).
 In particular, $\alpha_{x,y}(+\infty)$
and $\alpha_{y,x}(+\infty)$ belong to $\Lambda_{C'}$.
\end{lemma}

{\bf Proof:} Let $S_{x}$ and $S_{y}$ be the unit spheres of $E_{x}$ and $E_{y}$.
For any unit tangent vector $u\in T_{z}\tilde{X}$ at some point $z$,
we define $\theta_{u} \in \partial \tilde{X}$ by 
$\dot{\alpha}_{z,\theta_{u}}(0) = u$. 
By step 3, $\Lambda_{C'} = E_{x}(\infty) = E_{y}(\infty)$, therefore for every $u\in S_{x}$,
$\theta_{u}\in \Lambda_{C'}$ and 
there exist $v\in E_{y}$ such that $\theta_{u} = \theta_{v}$.
As $E_{y}$ is a vector space, $\theta_{-v}$ belongs to $\Lambda_{C'}$ therefore there exist
$w \in E_{x}$ such that $\theta_{w} = \theta_{-v}$.
The map $f: S_{x} \to S_{x}$ defined by
$f(u)=w$ is a continuous map.
The lemma then boils down to proving that there exist $u\in S_{x}$
such that $f(u)= -u$ because in that case, $x$, $y$ and $\theta_{u}$ are
on the same geodesic $\alpha_{x,\theta_{u}}$.

The following properties of $f$ are obvious. 

(i) For every $u\in S_{x}$, $f(u) \neq u$.

(ii) $f \circ f = Id$.

So $f$ is an involution of the sphere without fixed point and for any 
such map, we claim that there exist $u$ in the sphere such that $f(u)=-u$.
In order to prove the claim, we follow a very similar argument in 
\cite{sz}, theorem 1.
 We argue by contradiction.
Let us assume that for every $u\in S_{x}$, $f(u) \neq -u$.
The map $g: S_{x} \to S_{x}$ defined by
$g(u) = \frac{f(u) +u}{\|f(u)+u\|}$, is then well defined 
and continuous. 
Let us  remark that as for every $u\in S_{x}$, $f(u) \neq -u$, then $f$ is homotopic to the
Identity, and so is $g$. Moreover by (ii) we clearly have $g\circ f = g$, thus the map $g$
factorizes through $S_{x}/ G_{f}$ where $G_{f}$ is the group generated by the involution $f$.
By (i) $f$ has no fixed point thus $S_{x}/G_{f}$ is a manifold and 
the projection $p: S_{x} \to S_{x}/G_{f}$ is a degre 2 map. Therefore,
the induced endomorphism $g_{\ast}$ on $H_{n-1}(S_{x}, \mathbb{Z}_{2})$ is trivial,
which contradicts the fact that $g$ is homotopic to the Identity. $\Box$ 

\begin{cor}
Let $\mathcal{H}^{n-1} \subset \tilde{X}$ be an hypersurface with possibly
non empty boundary $\partial \mathcal{H}^{n-1}$, such that for any 
$x\in \mathcal{H}^{n-1}$, 
$Jac_{n-1}\tilde{F'}(x) = Jac_{E_{x}}\tilde{F'}(x) =1$, where
$E_{x}$ is the tangent space of $\mathcal{H}^{n-1}$ at $x$.
Let us consider $x\in \mathcal{H}^{n-1}$ 
such that $dist_{\tilde{X}}(x,\partial \mathcal{H}^{n-1})=r>0$.
Then, for any $x'\in \mathcal{H}^{n-1}$ with 
$dist_{\tilde{X}}(x,x') < r$, the geodesic $\alpha_{x,x'}$ joining $x$
and $x'$ is contained in $\mathcal{H}^{n-1}$. In particular, 
$\mathcal{H}^{n-1}$ is locally convex.
\end{cor}

{\bf Proof of the corollary:} Let us fix $\theta \in \Lambda_{C'}$ and consider
the vector field $\nabla B(y,\theta)$ in $\tilde{X}$. Let $x\in \mathcal{H}^{n-1}$. As 
$Jac_{E_{x}}\tilde{F'}(x) =1$, we have by step 3 $\Lambda_{C'} = E_{x}(\infty)$.
Then, for any $x\in \mathcal{H}^{n-1}$, $\nabla B(x,\theta)$ is tangent to 
$\mathcal{H}^{n-1}$, therefore the geodesic $\alpha_{x,\theta}$ satisfies
$\alpha_{x,\theta}(t) \in \mathcal{H}^{n-1}$ for all $t\in [0,r)$.
Let $x'\in \mathcal{H}^{n-1}$. By lemma 5.29, $\dot{\alpha}_{x,x'}(0) \in E_{x}$, therefore
$\alpha_{x,x'} =\alpha_{x,\theta}$ and $\alpha_{x,x'}(t) \in \mathcal{H}^{n-1}$ for
all $t\in [0,r)$. $\Box$

We now prove that $Z_{\infty}$ is a totally geodesic hypersurface in $\tilde{X}/C'$.
Let us recall that $Z_{\infty}$ is a manifold which is smooth except at a singular 
subset of codimension at least $7$. Let us consider a lift 
$\tilde{Z}_{\infty} \subset\tilde{X}$ of $Z_{\infty}$ and denote
$\tilde{Z}_{\infty}^{reg}$ (resp. $\tilde{Z}_{\infty}^{sing}$)
the set of regular (resp.) singular points of
$\tilde{Z}_{\infty}$.

\begin{lemma}
$\tilde{Z}_{\infty}$ is a totally geodesic hypersurface in $\tilde{X}$.
\end{lemma} 
 
{\bf Proof:} Let us consider a regular point $x\in \tilde{Z}_{\infty}^{reg}$. We shall show
that for every point $x'\in \tilde{Z}_{\infty}^{reg}$ the geodesic segment joining $x$ and $x'$
is contained in $\tilde{Z}_{\infty}$, and as the set of regular points is dense 
in $\tilde{Z}_{\infty}$ (as the complementary of
a subset of codimension at least $8$), this will show
that $\tilde{Z}_{\infty}$ is totally geodesic.

\noindent
We claim that there exist a sequence $y_{k} \in \tilde{Z}_{\infty}^{reg}$ such that
$\lim_{k\to \infty} y_{k} = x'$ and the geodesic segment joining $x$ and $y_{k}$
is contained in $\tilde{Z}_{\infty}$.

\noindent
The claim immediately implies that the geodesic segment joining $x$ and $x'$ 
is contained in $\tilde{Z}_{\infty}$.

Let us prove the claim.

For $y\in \tilde{Z}_{\infty}^{reg}$ we consider $\alpha_{x,y}$ the geodesic joining
$x$ and $y$ and define 

\begin{equation}
t_{y} = inf \{t >0 , \alpha_{x,y}(t)  \notin \tilde{Z}_{\infty} \}
\end{equation}  

As $x$ is a regular point,
by corollary 5.30, there exist $\epsilon >0$ such that $t_{y} > \epsilon$.

In order to prove the claim, we argue by contradiction. Let us assume
that there exist $r>0$ such that for any $y \in B_{\tilde{X}}(x',r)\cap \tilde{Z}_{\infty}^{reg}$,
$t_{y} < dist(x, y)$. By corollary 5.30 applied to $\tilde{Z}^{reg}_{\infty}$,
 we have $\alpha_{x,y}(t_{y}) \in \tilde{Z}_{\infty}^{sing}$.
 As the set of regular points is an open subset of $\tilde{Z}_{\infty}$,
if $r$ is small enough we have
$B_{\tilde{X}}(x',r) \cap \tilde{Z}_{\infty}^{reg}
= B_{\tilde{X}}(x',r) \cap \tilde{Z}_{\infty}$. We choose such an $r$ and 
we consider the set $S$ of all singular points contained in the union of 
all geodesic segments joining $x$ to a point $y\in B_{\tilde{X}}(x',r) \cap {Z}_{\infty}$.
Let us consider the map defined on $S$ by 

$$
p(y) = \alpha_{x,y}(\epsilon).
$$

As we already saw, for any $y\in B_{\tilde{X}}(x',r) \cap \tilde{Z}_{\infty}$,
we have $t_{y} > \epsilon$, therefore the map $p$ is distance decreasing and by assumption
$p$ is surjective onto an open subset of the sphere and $p(S)$ is homeomorphic to an open subset of $\mathbb {R}^{n-1}$,
therefore the Hausdorff dimension of $S$ is greater than or equal to $n-1$, which 
contradicts the fact that the singular set has codimension at least $8$ in $\tilde{Z}_{\infty}$.
$\Box$

The totally geodesic hypersurface $\tilde{Z}_{\infty}\subset \tilde{X}$ is preserved by $C$,
and $\tilde{Z}_{\infty} /C$ is of minimal volume in its homology class.

Let us prove that $\tilde{Z}_{\infty}$ is isometric to the hyperbolic space
$\mathbb{H}_{\mathbb{R}}^{n-1}$.

\begin{lemma}
$\tilde{Z}_{\infty}$ is isometric to the hyperbolic space $\mathbb{H}_{\mathbb{R}}^{n-1}$.
\end{lemma}

{\bf Proof :}

As $\tilde{Z}_{\infty}/C$ is of minimal volume in its homology class, we have by Proposition 5.1,
for all $x\in \tilde{Z}_{\infty}$, $Jac_{E_{x}}\tilde{F}(x) =1$ and 
$\tilde{F}(x) = x$, where $E_{x}$ is the tangent space of $\tilde{Z}_{\infty}$ at $x$.
Moreover, we saw in the proof of proposition 5.1 that

$$
H = \frac{1}{n-1}Id_{D\tilde{F}(x)(E_{x})} = \frac{1}{n-1} Id _{E_{x}},
$$

therefore we get from (4.11) and $\tilde{F}(x) = x$,
 that for all $u , v \in T_{x} \tilde{Z}_{\infty}$,

$$
\int_{\partial \tilde{X}}[ DdB_{(x,\theta)}(u,v) + DB_{(x,\theta)}(u) DB_{(x,\theta)}(v)]
 d\nu_{x}(\theta) 
$$

\begin{equation}
= \tilde{g}(u,v)
\end{equation}

where $\tilde{g}$ is the metric on $\tilde{X}$.
As $\tilde{Z}_{\infty}$ is totally geodesic, the relation (5.62) remains true with 
the Busemann function $B^{\tilde{Z}_{\infty}}$ of $\tilde{Z}_{\infty}$ instead
of the Busemann function $B$ of $\tilde{X}$:

$$
\int_{\partial \tilde{X}} [ DdB^{\tilde{Z}_{\infty}}_{(x,\theta)}(u,v) +
DB^{\tilde{Z}_{\infty}}_{(x,\theta)}(u)
DB^{\tilde{Z}_{\infty}}_{(x,\theta)}(v)
$$

\begin{equation}
= \tilde{g}(u,v).
\end{equation}

On the other hand, as $\tilde{Z}_{\infty}$ is totally geodesic, its sectional curvature
is less than or equal to $-1$, thus by Rauch comparison theorem, we have

\begin{equation}
DdB^{\tilde{Z}_{\infty}}(x,\theta) + DB^{\tilde{Z}_{\infty}}(x,\theta)
\otimes DB^{\tilde{Z}_{\infty}}(x,\theta)
\geq \tilde{g}_{\tilde{Z}_{\infty}}
\end{equation}

for all $\theta \in \partial \tilde{Z}_{\infty} = \Lambda_{C}$, where
$\tilde{g}_{\tilde{Z}_{\infty}}$ is the restriction of $\tilde{g}$ to 
$\tilde{Z}_{\infty}$.

As the support of the measure $\nu_{x}$ is $\partial \tilde{Z}_{\infty} = \Lambda_{C}$
(by convex cocompactness of $C$) and the Busemann function is continuous,
we get from (5.64) and (5.65) that 
for all $x\in \tilde{Z}_{\infty}$ and all $\theta \in \tilde{Z}_{\infty}$

$$
DdB^{\tilde{Z}_{\infty}}(x,\theta) +
DB^{\tilde{Z}_{\infty}}(x,\theta) \otimes DB^{\tilde{Z}_{\infty}}(x,\theta)
$$
 
\begin{equation}
= \tilde{g}_{\tilde{Z}_{\infty}}(x).
\end{equation}

and this last relation is characteristic of the hyperbolic space. $\Box$

\bigskip
{\bf Step 6: Conclusion}

\bigskip
So far we have shown that $C$ preserves a totally geodesic copy of the hyperbolic
space $\mathbb {H}_{\mathbb{R}}^{n-1} \subset \tilde{X}$
such that $\mathbb{H}^{n-1}_{\mathbb{R}} /C$ is compact.

Our goal now is to show that $Y =:\mathbb{H}^{n-1}_{\mathbb{R}}/C$ injects
diffeomorphically in $ X = \tilde{X}/\Gamma$ and separates $X$ in two 
connected components $R$ and $S$ such that $A = \pi_{1}(R)$ and $B = \pi_{1}(S)$.

In order to do this, we will consider the $\Gamma$ orbit of $\mathbb{H}^{n-1}_{\mathbb{R}}$
in $\tilde{X}$ and the two connected components $U$ and $V$ of
$\tilde{X} - \Gamma \mathbb{H}^{n-1}_{\mathbb{R}}$ which are adjacent to
$\mathbb{H}^{n-1}_{\mathbb{R}}$. The stabilizers $\bar{A}$, $\bar{B}$ and $\bar{C}$
of $U$, $V$ and $\mathbb{H}^{n-1}_{\mathbb{R}}$ contain respectively $A$, $B$ and $C$
and the hypersurface $\mathbb{H}^{n-1}_{\mathbb{R}} /C$ injects in $X = \tilde{X}/\Gamma$
and separates $X$ in two connected components $R$ and $S$ such that
$\pi_{1}(R)=\bar{A}$ and $\pi_{1}(S)=\bar{B}$.
We then show that $\bar{C} =C$, $\bar{A}=A$ and $\bar{B}=B$.

\bigskip
Let $\bar{C}$ be the stabilizer of $\mathbb H^{n-1}$, namely
$\bar{C} = \{\gamma \in \Gamma , \gamma \mathbb H^{n-1} =\mathbb H^{n-1}\}$.
We have $C\subset \bar{C}$ and as $\mathbb H^{n-1} /C$ is compact,
so is $\mathbb H^{n-1} /\Bar{C}$ and thus [$\bar{C} : C$] $< \infty$.

Let $p : \tilde{X}/C \rightarrow X = \tilde{X}/\Gamma$ 
and $\bar{p} : \tilde{X} /\bar{C} \rightarrow X= \tilde{X}/\Gamma$ 
the natural projections.
We now show that the restriction of $p$ to $\mathbb{H}^{n-1}/C$ is an embedding,
thus $Y:=p(\mathbb{H}^{n-1}/C)$ is a compact
 totally geodesic hypersurface of $X$.

In the section 2, we constructed a $C$-invariant hypersurface 
$\tilde{Z} \subset \tilde{X}$ such that $Z = \tilde{Z}/C \subset \tilde{X}/C$
is compact.
The hypersurface is defined as 
$\tilde{Z} = \tilde{f}^{-1}(t_{0})$
where $\tilde{f} : \tilde{X} \rightarrow T$ is an
equivariant map onto the Bass-Serre tree associated to the amalgamation 
$A\ast_{C}B$ and $t_{0}$ belongs to that edge of $T$
which is fixed by $C$.

Let us first show two lemmas.

\begin{lemma}
The restriction of $p$ to $\tilde{Z}/C$ is an embedding 
into $X=\tilde{X}/\Gamma$.
\end{lemma}

{\bf Proof :} Let $\gamma \in \Gamma$, $z$, $z'$ in $\tilde{Z}$ such that
$z'=\gamma z$. By equivariance, 

$$
\tilde{f}(\gamma z) = \gamma \tilde{f}(z) = \gamma t_{0} =\tilde{f}(z') = t_{0},
$$
thus $\gamma \in C$.$\Box$

\begin{lemma}
The restriction of $\bar{p}$ to $\mathbb H^{n-1}/\bar{C}$ is an embedding 
into $X=\tilde{X}/\Gamma$.
\end{lemma}

{\bf Proof :} Let us assume that there is a $\gamma \in \Gamma - \bar{C}$
such that $\gamma \mathbb H^{n-1} \cap \mathbb H^{n-1} \neq \emptyset$
and choose an $x\in \gamma \mathbb H^{n-1} \cap \mathbb H^{n-1}$.
As $\gamma \notin \bar{C}$, there exist 
$u \in T_{x} \gamma \mathbb H^{n-1} - T_{x}\mathbb H^{n-1}$.
We consider $c_{u}$ the geodesic ray such that $\dot{c}_{u}(0) = u$.
We know that $\tilde{Z}$ is contained in an 
$\epsilon$-neighbourhood  $\mathcal {U} _{\epsilon} \mathbb H^{n-1}$ 
 of $\mathbb H^{n-1}$. 
The $\epsilon$-neighbourhood  $\mathcal {U} _{\epsilon} \mathbb H^{n-1}$ 
 of $\mathbb H^{n-1}$ separates $\tilde{X}$ in two 
connected components $U$ and $V$ and for $t>0$ large enough,
we have, say,  $c_{u}(t)\in U$ and $c_{u}(-t)\in V$. 

Let $\tilde{Z'}$ be the connected component of $\tilde{Z}$ that we constructed
at the end of section 2, whose stabilizer 
(or an index two subgroup of it) $C'$ is such that $\tilde{Z'}/C'$
separates $\tilde{X}/C'$ in two unbounded connected components
$U'/C'$ and $V'/C'$ where $U'$ and $V'$ are 
the two connected components of $\tilde{X} - \tilde{Z'}$.

We claim that $U \subset U'$ and $V \subset V'$ or the other way around.
Indeed if not, $U$ and $V$ would be both contained in, say, $U'$. But in that case,
$V'$ would be contained in $\mathcal{U}_{\epsilon} H^{n-1}$ and therefore
$V'/C'$ would be bounded, which is a contradiction. 
 
As $\gamma \tilde{Z}$ lies in the $\epsilon$ neighborhood of
 $\gamma \mathbb H^{n-1}$,
there exist sequences $z_{k}$, $z'_{k}$ in $\gamma \tilde{Z}$
such that $dist(z_{k},c_{u}(k)) \leq \epsilon$
and  $dist(z'_{k},c_{u}(-k)) \leq \epsilon$.
By proposition 5.13 and lemma 5.23, 
$C'$ also acts
cocompactly on $\mathbb H^{n-1}$, thus $C'$ is of finite index in $C$,
and therefore there are finitely many connected components of $\tilde{Z}$
and the same holds for $\gamma \tilde{Z}$.
We thus can assume that the $z_{k}$'s and $z'_{k}$'s belong to a single
connected component of $\gamma \tilde{Z}$.
Let us consider a continuous path $\alpha \subset \gamma \tilde{Z}$
 joining $z_{k}$ and $z'_{k}$.

By construction the distance between $c_{u}(k)$ [resp. $c_{u}(-k)$] and 
$\mathbb H^{n-1}$ tends to infinity and thus, for
$k$ large enough, $z_{k}\in U$ and $z'_{k}\in V$
or the other way around. By the claim, we then have $z_{k}\in U'$ and $z'_{k}\in V'$,
therefore the path $\alpha$ has to cross $\tilde{Z'}$ which contradicts the lemma (5.29)
and ends the proof of the lemma 5.30. $\Box$

As we already saw, $\tilde{Z}$ has finitely many connected components,
and so does $\tilde{X} - \tilde{Z}$. 
Let us write $\{W_{j}\}_{j=1,..,m}$ the connected components of $\tilde{X} - \tilde{Z}$.
As $C$ acts cocompactly on $\tilde{Z}$ and $\mathbb H^{n-1}$ there exist 
$\epsilon >0$ such that $\mathbb H^{n-1} \subset \mathcal {U}_{\epsilon} \tilde{Z}$
and $\tilde{Z} \subset \mathcal{U}_{\epsilon}\mathbb H^{n-1}$. Moreover 
$\mathcal{U}_{\epsilon} \mathbb H^{n-1}$ separates $\tilde{X}$ in two
connected components $U$ and $V$. 

\begin{lemma}
Let us consider $\epsilon$ such that
$\tilde{Z} \subset \mathcal{U}_{\epsilon}\mathbb H^{n-1}$ and
$U$ and $V$ the two connected components of 
$\tilde{X} - \mathcal{U}_{\epsilon}\mathbb H^{n-1}$.
There are two distinct connected components 
$W_{1}$ and $W_{2}$ of $\tilde{X} - \tilde{Z}$ 
such that $U\subset W_{1}$ and $V\subset W_{2}$.
Moreover, $\tilde{f}(W_{1}) \subset \tilde{T}_{1}$
and $\tilde{f}(W_{2}) \subset \tilde{T}_{2}$, where
$\tilde{T}_{1}$ and $\tilde{T}_{2}$ are the two connected components
of $\tilde{T} -\{t_{0}\}$.
\end{lemma}

{\bf Proof :} We argue by contradiction. Let us assume 
that $U$ and $V$ are contained in the same connected component $W_{1}$
of $\tilde{X} - \tilde{Z}$. Then, all other components $W_{j}$, $j\neq 1$, satisfy
$W_{j} \subset \mathcal {U}_{\epsilon} \mathbb H^{n-1} \subset \mathcal {U}_{2\epsilon}\tilde{Z}$.
Therefore, as $C$ acts cocompactly on $\mathcal {U}_{2\epsilon}\tilde{Z}$, 
there exist a constant $D$ such that for any $j \neq 1$,
$max _{w\in W_{j}} dist _{\tilde{T}}(\tilde{f}(w),t_{0}) \leq D$.
Thus, $\tilde{f}(W_{1})$ is contained in one connected component of 
$\tilde{T}-\{t_{0}\}$ and 
$\tilde{f}(\cup _{j\neq 1} W_{j})$, contained in
the ball $B_{\tilde{T}}(t_{0}, D)$
of $\tilde{T}$ of radius $D$ centered at $t_{0}$, is bounded. This is clearly
impossible because $\tilde{T}-\{t_{0}\}$ has two unbounded connected components
and $\tilde{f}$ is onto. $\Box$ 

Let us denote $\mathcal {A} = A \mathbb H^{n-1}$ the $A$-orbit of 
the $C$-invariant totally geodesic copy  of the real hyperbolic space $\mathbb H^{n-1}$,
and $\bar{A}$ the stabilizer of $\mathcal {A}$, ie. 
$\bar {A}= \{\gamma \in \Gamma , \gamma \mathcal{A} = \mathcal{A} \}$.
We define in a similar way $\mathcal{B} = B \mathbb H^{n-1}$ and 
$\bar{B} = \{\gamma \in \Gamma , \gamma \mathcal {B} = \mathcal {B}\}$.

Let us recall that $\bar{C}$ is the stabilizer of $\mathbb H^{n-1}$ in $\Gamma$.
We now prove the following

\begin{lemma}
We have $\bar{A} = A \bar{C}$ and $\bar{B} = B \bar{C}$. Moreover,
$\bar{A}$ and $\bar{B}$ are charactrized by 
$\bar{A} = \{ \gamma \in \Gamma , \gamma \mathbb H^{n-1} \in \mathcal {A} \}$ and
$\bar{B} = \{ \gamma \in \Gamma , \gamma \mathbb H^{n-1} \in \mathcal {B} \}$.
\end{lemma}

{\bf Proof :} Let $\gamma '\in \bar{A}$, then 
$\gamma '\mathbb H^{n-1} \in \mathcal{A}$ and 
thus there exist $\gamma  \in A$ such that 
$\gamma '\mathbb H^{n-1} = \gamma  \mathbb H^{n-1}$, therefore $\gamma ^{-1}\gamma '\in \bar{C}$,
which proves the first part of the lemma.

Let us prove the second part of the lemma.

Let $\gamma ' \in \Gamma$ be such that $\gamma '\mathbb H^{n-1} \in \mathcal{A}$.
Then there exist $\gamma \in A$ such that 
$\gamma ' \mathbb H^{n-1} = \gamma \mathbb H^{n-1}$, thus 
$\gamma^{-1} \gamma ' \in \bar{C}$ and therefore
$\gamma ' \in \gamma \bar{C} \subset A\bar{C} = \bar{A}$.
This proves one inclusion, the other inclusion being obvious. $\Box$

For each $\gamma \in \Gamma$, $\gamma \mathbb H^{n-1}$ separates
$\tilde{X}$ in two connected components $U_{\gamma}$ and $V_{\gamma}$.

Let us now prove the following lemma.

\begin{lemma}

\noindent
(i) Let $\gamma \in A$, [resp. $\gamma \in B$]. Then, we have
 $\mathcal{A} - \{ \gamma \mathbb H^{n-1}\} \subset U_{\gamma}$
 or $\mathcal{A} -\{ \gamma \mathbb H^{n-1} \} \subset  V_{\gamma}$,
[resp. 
$\mathcal{B} - \{ \gamma \mathbb H^{n-1}\} \subset U_{\gamma}$
 or $\mathcal{B} -\{ \gamma \mathbb H^{n-1} \} \subset  V_{\gamma}$.]

(ii) Let $\gamma$ be an element of $\Gamma - \bar{A}$, [resp. $\Gamma - \bar{B}$].
Then $\mathcal{A} \subset U_{\gamma}$ or 
$\mathcal{A} \subset V_{\gamma}$,
[resp. $\mathcal{B} \subset U_{\gamma}$ or
$\mathcal{B} \subset V_{\gamma}$].
\end{lemma}

{\bf Proof :} (i) We argue by contradiction. Let us consider 
$\gamma \mathbb H^{n-1}$, $\gamma ' \mathbb H^{n-1}$ and
$\gamma '' \mathbb H^{n-1}$ three distinct elements in $\mathcal {A}$
such that $\gamma' \mathbb H^{n-1} \subset U_{\gamma}$ and
$\gamma '' \mathbb H^{n-1} \subset V_{\gamma}$. 
By equivariance we can assume $\gamma$ is the identity.
Let us recall that $U$ and $V$ are the two connected components of
$\tilde{X} - \mathcal{U}_{\epsilon}\mathbb H^{n-1}$.

We then have $\gamma ' \mathbb H^{n-1} \cap U \neq \emptyset$
and $\gamma '' \mathbb H^{n-1} \cap V \neq \emptyset$, which implies
$\gamma ' \tilde{Z} \cap U \neq \emptyset$ and 
$\gamma '' \tilde{Z} \cap V \neq \emptyset$.

By lemma (5.31), $U \subset W_{1}$ and $V\subset W_{2}$
where $W_{1}$ and $W_{2}$ are two connected components of 
$\tilde{X}-\tilde{Z}$ and 
$\tilde{f}(U)\subset \tilde{T}_{1}$ 
and $\tilde{f}(V) \subset \tilde{T}_{2}$,
therefore, $\tilde{f}(U)$ contains $\gamma ' t_{0} \in \tilde{T}_{1}$
and $\tilde{f}(V)$ contains $\gamma '' t_{0}\in \tilde{T}_{2}$. 
This is impossible because for all elements $\gamma '$ 
and $\gamma ''$ in $A$, 
$\gamma ' t_{0}$ and $\gamma '' t_{0}$ belong to the same 
connected component of $\tilde{T} - \{t_{0} \}$. 

(ii) Let us consider $\gamma \in \Gamma - \bar{A}$.
We argue by contradiction.
Let us assume there exist $\gamma$, $\gamma'$ in $A$ such that

$$
\gamma'\mathbb H^{n-1} \subset U_{\gamma} 
$$
\begin{equation}
\gamma'' \mathbb H^{n-1} \subset V_{\gamma}
\end{equation}

Let $\epsilon >0$ such that 
$\tilde{Z} \subset \mathcal{U}_{\epsilon}\mathbb H^{n-1}$
and $U$ and $V$ the connected component of $\tilde{X} - \mathcal{U}_{\epsilon}\mathbb H^{n-1}$.
By lemma (5.31) we have 
$\tilde{f}(\gamma U) \subset \gamma \tilde{T}_{1}$
and
$\tilde{f}(\gamma V) \subset \gamma \tilde{T}_{2}$,
where 
$\gamma \tilde{T}_{1}$ and $\gamma \tilde{T}_{2}$
are the two connected components of
$\tilde{T}-\{\gamma t_{0} \}$.
By assumption (5.62), we have 

$\gamma ' \mathbb H^{n-1} \cap \gamma U \neq \emptyset$
and $\gamma '' \mathbb H^{n-1} \cap \gamma V \neq \emptyset$, which implies
$\gamma ' \tilde{Z} \cap \gamma U \neq \emptyset$ and 
$\gamma '' \tilde{Z} \cap \gamma V \neq \emptyset$,
therefore 
$\gamma' t_{0} \in \gamma \tilde{T}_{1}$ and
$\gamma'' t_{0} \in \gamma \tilde{T}_{2}$,
which is impossible because in the tree $\tilde{T}$,
the points $\gamma' t_{0}$ and $\gamma'' t_{0}$ belong to two adjacent edges.

$\Box$ 

By lemma (5.33) (i), for every 
$\gamma$ in $\bar{A}$, [resp. $\bar{B}$ ], we can define $U_{\gamma}$ as the connected component 
of $\tilde{X}-\mathbb H^{n-1}$ which contains all $\gamma ' \mathbb H^{n-1}$
for all $\gamma'$ in $\bar{A}$, [ resp. $\bar{B}$ ], 
 and $\gamma' \mathbb H^{n-1} \neq \gamma \mathbb H^{n-1}$.

Let us define
 
\begin{equation}
U_{A}:= \cap _{\gamma \in A} U_{\gamma}. 
\end{equation}

By  definition, $U_{A}$ [resp. $U_{B}$ ]is a convex set in $\tilde{X}$ 
whose boundary is the collection $\mathcal{A}$, [resp. $\mathcal{B}$ ],
of $\gamma \mathbb H^{n-1}$,
$\gamma$ in $\bar{A}$ [resp. $\bar{B}$ ] and by
lemma (5.33) (ii), $U_{A}$ and $U_{B}$ are two disjoint connected components of 
$\tilde{X} - \Gamma \mathbb H^{n-1}$.

In fact, $U_{A}$ [resp. $U_{B}$], is the convex hull of $\mathcal{A}$, [resp. $\mathcal{B}$],
and $\bar{A}$, [resp. $\bar{B}$], is the stabilizer of $U_{A}$, [resp. $U_{B}$]. 

\begin{lemma}
The closures of 
$U_{A}$ and $U_{B}$ intersect along $\mathbb H^{n-1}$ and
$\bar{A} \cap \bar{B} = \bar{C}$. Moreover, no element
of $\Gamma$ sends $U_{A}$ on $U_{B}$ nor the other way around.
\end{lemma}

{\bf Proof :}  The convex set $U_{A}$ is the intersection of open half spaces
$U_{\gamma}$, $\gamma \in \bar{A}$, and is delimited by the disjoint union of hyperplanes
$\gamma \mathbb H^{n-1}$, for some $\gamma \in \bar{A}$.
The same is true for $U_{B}$ and as $U_{A} \cap U_{B} = \emptyset$,
the closures of $U_{A}$ and $U_{B}$ can intersect 
only along one of the connected components of their boundaries, thus 
along $\mathbb H^{n-1}$ which is obviously in both closures.
This proves the first part of the lemma, let us prove the 
second part.
By lemma 5.32, $\bar{C} \subset \bar{A} \cap \bar{B}$.
Conversely, let us take $\gamma \in \bar{A} \cap \bar{B}$, then $\gamma$ preserves
the closures of $U_{A}$ and $U_{B}$, thus it preserves their intersection $\mathbb H^{n-1}$,
and therefore $\gamma \in \bar{C}$.

Let us prove the last part of the lemma.
Let $\gamma$ be an element such that 
$\gamma U_{A} = U_{B}$. 
As $\mathbb H^{n-1}$ is one component of the boundary $\mathcal{B}$ of $U_{B}$,
there exist one component $\gamma' \mathbb H^{n-1} \in \mathcal{A}$, 
$\gamma'$ being in $\bar{A}$, such that 
$\gamma (\gamma' \mathbb H^{n-1}) = \mathbb H^{n-1}$.
Therefore, $\gamma \gamma' \in \bar{C}$, thus $\gamma \in \bar{A}$.
The same argument yields $\gamma^{-1} \in \bar{B}$, so 
$\gamma \in \bar{A} \cap \bar{B} = \bar{C}$ and 
$\gamma$ preserves $U_{A}$ and $U_{B}$, which contradicts
our choice of $\gamma$.

 $\Box$

The $\Gamma$-orbit of the closure of $U_{A} \cup U_{B}$ covers
$\tilde{X}$.
Let us construct a tree $\bar{T}$ embedded in $\tilde{X}$ in the following way:
the set of vertices is the set of the connected components of 
$\tilde{X} - \Gamma \mathbb H^{n-1}$ and two vertices
are joined by an edge if the boundaries 
of their corresponding connected components 
intersect non trivially in $\tilde{X}$.
By construction $\Gamma$ acts on $\bar{T}$, the stabilizers of the vertices 
$a$ and $b$ corresponding to $U_{A}$ and $U_{B}$ 
are $\bar{A}$ and $\bar{B}$, the stabilizer of the edge between $a$ and $b$
is $\bar{C}$ and a fundamental domain for this action is the segment
joining $a$ and $b$. By \cite{se}, I, 4, Theorem 6, the group 
$\Gamma$ is the amalgamated product of $\bar{A}$ and $\bar{B}$ over $\bar{C}$.

We now claim that $\bar{A}=A$, $\bar{B}=B$ and $\bar{C}=C$.

As $A$, $B$ and $C$ are subgroups of $\bar{A}$, $\bar{B}$, and $\bar{C}$,
the corresponding Mayer-Vietoris sequences of $A \ast_{C} B$ and 
$\bar{A}\ast _{\bar{C}} \bar{B}$ are related by the following commutative
diagram

$$
\begin{matrix}
& H_{n}(A,\mathbb R) \oplus H_{n}(B,\mathbb R) &
 \longrightarrow & H_{n}(\Gamma,\mathbb R)&
  \longrightarrow & H_{n-1}(C,\mathbb R) \\
&\downarrow & & \downarrow & & \downarrow \\
& H_{n}(\bar{A},\mathbb R) \oplus H_{n}(\bar{B},\mathbb R) &
 \longrightarrow & H_{n}(\Gamma,\mathbb R)&
  \longrightarrow & H_{n-1}(\bar{C},\mathbb R). \\
\end{matrix}
$$

We know that the index  $[\bar{C}:C]$ is finite, and 
 by the lemma 5.32, the indices of 
$A$, $B$, and $C$ in  $\bar{A}$, $\bar{B}$ and $\bar{C}$ are finite and equal.
On the other hand, the indices $[\Gamma : A]$ and $[\Gamma : B]$ are infinite
by assumption, thus the previous diagram becomes

$$
\begin{matrix}
&0& \longrightarrow & H_{n}(\Gamma,\mathbb R)&
  \longrightarrow & H_{n-1}(C,\mathbb R) \\
& & & \downarrow & & \downarrow \\
&0& \longrightarrow & H_{n}(\Gamma,\mathbb R)&
  \longrightarrow & H_{n-1}(\bar{C},\mathbb R). \\
\end{matrix}
$$

Moreover, the map $H_{n}(\Gamma, \mathbb{R}) \to H_{n-1}(\bar{C},\mathbb{R})$
is bijective. Namely, the injectivity comes from the above diagram and the surjectivity 
from the fact that that the hypersurface $\mathbb{H}^{n-1}/\bar{C}$ bounds in 
$\mathbb{H}^{n}/A$ and  $\mathbb{H}^{n}/B$ so that the map 
$H_{n-1}(C,\mathbb{R}) \to H_{n-1}(A,\mathbb{R}) \oplus H_{n-1}(B,\mathbb{R}) $ is trivial.
Therefore the index $[\bar{C} : C] =1$ and we get 
 $\bar{A}=A$, $\bar{B}=B$ and $\bar{C}=C$.
$\Box$

\section{Proof of the theorems 1.5 and 1.6.}

The proof of theorem 1.5 is exactly the same as the proof of theorem 1.2.
The actions of $\Gamma$ on $\tilde(X)$ and $T$ give rise to a continuous $\Gamma$-equivariant
 map $\tilde{f} : \tilde{X} \to T$. Like in section 2, we build an hypersurface
$\tilde{f}^{-1} (t_{0})$ where $t_{0}$ is a regular value of $\tilde{f}$ belonging the interior
of an edge. As the edge separates the tree in two unbounded components, the section 2 applies
and we get a subgroup $C'$ of $C$, and an hypersurface $\tilde{Z}' \subset \tilde{X}/C'$
which is essential.
Now, if the action of $\Gamma$ is minimal, every edge separates $T$ in two unbounded components.
$\Box$

\section{Appendix}

The goal of this section is to give a proof of lemma 5.9.
This lemma is contained in lemma 2.1, 5.1 and 5.2 of \cite{bk},
but our situation being not exactly the same,  we reproduce it down here
for sake of completness.

Let us restate the lemma 5.9.

\begin{lemma}
Let $\mathcal{L} \subset \partial \tilde{X}$ be a closed $C'$-invariant subset
and $\theta_{0} \in \mathcal{L}$. We assume that there exist
a sequence of positive
real numbers $\lambda_{k} \to\infty $
such that the sequence of pointed metric spaces 
$(\mathcal{L},\lambda_{k}d,\theta_{0})$ converges in the pointed Gromov-Hausdorff topology to    
$(S,\delta,0)$ where $(S,\delta,0)$ is a weak tangent of
 $(\partial \tilde{X},d)$.
We also assume that there exist positive constants $C$ and $\delta$, a sequence of 
points 
$\theta_{0}^{k}=\theta_{0}, \theta_{1}^{k}, \theta_{2}^{k} \in \mathcal{L}$
and a sequence of elements $\gamma_{k} \in C'$ such that 
$C^{-1} \lambda_{k}^{-1}\leq d(\theta_{i}^{k}, \theta_{j}^{k}) \leq  C \lambda_{k}^{-1}$
and $d(\gamma_{k}\theta_{i}^{k},\gamma_{k} \theta_{j}^{k}) \geq \delta$
for all $0\leq i\neq j \leq 2$. Then, $\mathcal{L}$ is homeomorphic to
the one point compactification $\hat{S}$ of $S$.
In particular $\mathcal{L}$ is homeomorphic to  $\partial \tilde{X}$.
\end{lemma} 

We first give a definition of pointed Hausdorff-Gromov convergence
which is equivalent to the definition 5.7.
We follow \cite{bk}, paragraph 4.

A sequence of metric spaces $(Z_{k}, d_{k}, z_{k})$ converges to the metric
space $(S, \delta, 0)$ if for every $R >0$, and every $\epsilon >0$, there exist an integer $N$,
a subset $D\subset B_{S}(0,R)$, subsets $D_{k} \subset B_{Z_{k}}(z_{k}, R)$ and 
bijections $f_{k} : D_{k} \rightarrow D$ such that for $k\geq N$,

(i) $f_{k} (z_{k}) = 0$,

(ii) the set $D$ is $\epsilon$-dense in $B_{S}(0,R)$, and the sets $D_{k}$

are $\epsilon$-dense in $B_{Z_{k}}(z_{k},R)$,

(iii) $|d_{Z_{k}}(x,y) - d_{Z}( f_{k}(x),f_{k}(y)) | < \epsilon$,

where $x$, $y$ belong to $D_{k}$.

Let us describe now the lemmas 2.1 and 5.1 following  \cite{bk}.

For a metric space $(Z,d)$ the cross ratio of four points $\{z_{i}\}$, $i=1,...4$, is the 
quantity


\begin{equation}
[z_{1},z_{2},z_{3},z_{4}] := \frac{d(z_{1},z_{3}) d(z_{2},z_{4})}{d(z_{1},z_{4}) d(z_{2},z_{3})}
\end{equation}

Given two metric spaces $X$ and $Y$, 
an homeomorphism $\eta : [0, \infty) \rightarrow [0, \infty)$,
 and an injective map $f : X \rightarrow Y$, we say that $f$ is an $\eta$-quasi-M\"{o}bius
map if for any four points $\{x_{i}\}$, $i=1,..,4$, in $X$,
we have 

\begin{equation}
[f(x_{1}),f(x_{2}),f(x_{3}),f(x_{4})] \leq \eta ( [x_{1},x_{2},x_{3},x_{4}]). 
\end{equation}

For example, any discrete cocompact group of isometries of $\tilde{X}$,
where $\tilde{X}$ is a Cartan-Hadamard manifold with sectional curvature
$K \leq -1$, is acting on the ideal boundary $(\partial \tilde{X}, d)$
endowed with the Gromov distance by $\eta$-quasi-M\"{o}bius transformations 
for some $\eta$.

\begin{lemma} 
[\cite{bk}, Lemma 2.1 ]
Let $(X,d_{X})$ and $(Y,d_{Y})$ be two compact metric spaces,
and for any integer $k$, $g_{k} : \tilde{D}_{k} \rightarrow Y$ 
an $\eta$-quasi-M\"{o}bius map defined on a subset $\tilde{D}_{k}$ of $X$.
We assume that the Hausdorff distance between $\tilde{D}_{k}$ and $X$ satisfies

$$
lim_{k\to \infty} dist_{H}(\tilde{D}_{k},X) = 0
$$
and that for any integer $k$, there exist points $(x_{1}^{k}, x_{2}^{k}, x_{3}^{k})$
in $D_{k}$ and  $(y_{1}^{k}, y_{2}^{k}, y_{3}^{k})$ in $Y$, such that 
$g_{k}(x_{i}^{k}) = y_{i}^{k}$ for $i\in\{1,2,3\}$, 
$d_{X}(x_{i}^{k},x_{j}^{k}) \geq \delta$ and
$d_{Y}(y_{i}^{k},y_{j}^{k}) \geq \delta$ for $i,j \in \{1,2,3\}, i\neq j$,
where $\delta$ is independant of $k$.
Then a subsequence of $g_{k}$ converges uniformly to 
a quasi-M\"{o}bius map $f :X \rightarrow Y$, ie. 
$lim_{k_{j}\to \infty}dist _{H}(g_{k_{j}}, f|_{\tilde{D}_{k_{j}}}) = 0$.
If in addition, we suppose that
$$
lim_{k\to \infty} dist_{H}(g_{k}(\tilde{D}_{k}),Y) = 0, 
$$
then the sequence $\{g_{k_{j}} \}$ converges uniformly to a quasi-M\"{o}bius 
homeomorphism $f: X \rightarrow Y$.
\end{lemma}

Before stating the second lemma, let us define a metric space $Z$ to be uniformly perfect
if there exist a constant $\lambda \geq 1$ such that for every $z\in Z$ and
$0 < R < diam Z$, we have 
$\bar{B}(z,R) - B(z,\frac{R}{\lambda})\neq \emptyset$. 

\begin{lemma}
[\cite{bk}, lemma 5.1 ]
Let $Z$ be a compact uniformly perfect metric space and $G$ an $\eta$-quasi-M\"{o}bius
action on $Z$.
Suppose that for each integer $k$ we are given a set $D_{k}$ in a ball
$B_{k} = B(z,R_{k}) \subset Z$ that is $(\epsilon_{k} R_{k})$-dense in $B_{k}$,
where $\epsilon_{k} >0$, distinct points 
$x_{1}^{k}, x_{2}^{k}, x_{3}^{k} \in B(z,\lambda_{k}R_{k})$, where $\lambda_{k} > 0$,
with 
$$
d_{Z}(x_{i}^{k},x_{j}^{k}) \geq \delta _{k} R_{k}
$$
for $i, j \in \{1, 2, 3\}, i \neq j$,
where $\delta_{k} >0$, and groups elements $\gamma_{k} \in G$ such that
for $y_{i}^{k} : = \gamma _{k} (x_{i}^{k})$ we have,
$$
d_{Z}( y_{i}^{k}, y_{j}^{k}) \geq \delta'
$$
for $i , j \in \{1, 2, 3 \}, i\neq j$, where $\delta'$ is independant of $k$.
Let $D'_{k} = \gamma _{k} (D_{k})$, and suppose that $\lambda_{k} \to 0$ when
$k\to \infty$, and the sequence $\frac{\epsilon_{k}}{\delta_{k}^{2}}$ is bounded.
Then $lim _{k\to \infty} dist_{H}(D'_{k}, Z) = 0$.
\end{lemma}

Let us go back to the proof of lemma 6.1.
By definition of convergence, there exist a subsequence of $\{\lambda_{k} \}$,
which we still denote by $\{\lambda_{k} \}$, subsets
$\tilde{D}_{k} \subset B_{S}(0, k)$,
$D_{k} \subset B_{\lambda_{k} \mathcal{L}}(\theta_{0}, k)$,
where $\tilde{D}_{k}$ and $D_{k}$ are minimal $1/k$-dense subsets of $B_{S}(0, k)$ and
$B_{(\mathcal{L},\lambda _{k} d)}(\theta _{0},k)$, and bijections
$f_{k} : \tilde{D}_{k} \rightarrow D_{k}$
such that for all $x, y \in \tilde{D}_{k}$,

\begin{equation}
\frac{1}{2} \delta(x,y) \leq \lambda_{k} d(f_{k}(x), f_{k}(y))
\leq 2 \delta (x,y),
\end{equation}

cf. \cite{bk}, (5.4).

We can suppose that the points $\theta_{0}^{k} : = \theta_{0}$, $\theta_{1}^{k}$, 
and $\theta_{2}^{k}$ in lemma 6.1 belong to the set $D_{k}$.
By assumption there exist elements $\gamma_{k} \in C'$ and a constant
 $\delta$ such that 

\begin{equation}
d(\gamma_{k} \theta_{i}^{k}, \gamma_{k} \theta_{j}^{k}) \geq \delta
\end{equation}  

for all $i, j \in \{0, 1, 2 \}$.

The lemma 6.1 is a direct consequence of the lemma 6.2 applied
to  $(X,d_{X}) = (\hat{S},\hat{\delta})$ and
$(Y,d_{Y}) = (\mathcal{L},d)$
 and to the sequence of maps
$g_{k} := \gamma _{k}\circ f_{k}$,
where $\hat{S}$ is the one point compactification of $S$ and
$\hat{\delta}$ the distance on $\hat{S}$ associated to $\delta$,
cf. \cite{bk} Lemma 2.2.

Let us denote $x_{0}^{k}, x_{1}^{k}, x_{2}^{k}$ be the points in $S$
such that $f_{k}(x_{i}^{k}) = \theta_{i}^{k}$, for $i \in \{0, 1, 2 \}$.

Let us check that the assumptions of lemma 6.2 are verified.

The fact that $lim_{k\to \infty} dist_{H}( \tilde{D}_{k}, \hat{S}) = 0$
comes the same way as in \cite{bk}, (5.5).

By (6.3), we have, 
$\delta( x_{i}^{k}, x_{j}^{k}) \geq \frac{\lambda_{k}}{2} 
d(\theta_{i}^{k}, \theta_{j}^{k})$
and by assumption we then get

\begin{equation}
\delta( x_{i}^{k}, x_{j}^{k}) \geq \frac{1}{2C}.
\end{equation}

We then get the separation assumption on triples of points
by choosing $\delta := inf \{ D, \frac{1}{2C}\}$.

It remains to check the assumption on
 $g_{k} (\tilde{D}_{k}) = \gamma_{k} \circ f_{k} (\tilde{D}_{k}) = \gamma_{k} (D_{k})$,
namely, 

\begin{equation}
 lim _{k\to \infty } dist_{H} (\gamma_{k} (D_{k}), \Lambda_{C'}) = 0.
\end{equation}

In order to prove the property (6.6), we want to apply the 
lemma 6.3, but as the set $(\mathcal{L},d)$ is 
a priori not uniformly perfect, we shall replace 
the uniform perfectness  by the fact that 
$(\mathcal{L}, \lambda_{k} d, \theta_{0})$
converges to a space $(S, \delta_{0}, 0)$, 
which is uniformly perfect, cf. ().

We will show the 

\begin{lemma} We consider the subsets
$\tilde{D}_{k} \subset B_{S}(0, k)$ and
$D_{k} \subset B_{\lambda_{k} \mathcal{L}}(\theta_{0}, k)$,
where $\tilde{D}_{k}$ and $D_{k}$ are $1/k$-dense subsets of $B_{S}(0, k)$ and
$B_{(\mathcal{L},\lambda _{k} d)}(\theta _{0},k)$, and the bijections
$f_{k} : \tilde{D}_{k} \rightarrow D_{k}$ coming from the convergence of 
 the sequence of pointed metric spaces 
$(\mathcal{L},\lambda_{k}d,\theta_{0})$  to    
$(S,\delta,0)$ where $(S,\delta,0)$ is a weak tangent of 
$(\partial \tilde{X},d)$.
We also assume that there exist positive constants $C$ and $\delta$, a sequence of 
points $\theta_{1}^{k}, \theta_{2}^{k} \in \Lambda_{C'}$
and a sequence of elements $\gamma_{k} \in C'$ such that 
$C^{-1} \lambda_{k}^{-1}\leq d(\theta_{i}^{k}, \theta_{j}^{k}) \leq  C \lambda_{k}^{-1}$
and $d(\gamma_{k}\theta_{i}^{k},\gamma_{k} \theta_{j}^{k}) \geq \delta$
for all $0\leq i\neq j \leq 2$. Then, the Hausdorf distance $dist_{H}(\gamma_{k} D_{k}, \mathcal{L})$
tends to $0$ as $k$ tends to infinity.
\end{lemma}

{\bf Proof :} The proof is word by word the same as the proof of 
lemma 6.3, ie. lemma 5.1 (i)
of \cite{bk} with a difference in case 2).

We have 
$B_{\lambda_{k} \mathcal{L}}(\theta_{0}, k) = B_{\mathcal{L}}(\theta_{0}, \frac{k}{\lambda_{k}})$
and $D_{k} \subset B_{\lambda_{k}\mathcal{L}}(\theta_{0}, k)$ an 
$\frac{1}{k}$-dense subset, for the metric $\lambda_{k} d$.
In term of the distance $d$, the set 
$D_{k}$ is $(\epsilon_{k} R_{k})$-dense in 
$B_{\mathcal{L}}(\theta_{0}, R_{k})$,
where $R_{k} := \frac{k}{\lambda_{k}}$ 
and $\epsilon_{k} := \frac{1}{k^{2}}$.
By assumption, the points $\theta_{0}^{k}=\theta_{0}, \theta_{1}^{k}, \theta_{2}^{k}$
belong to $B_{\mathcal{L}}(\theta_{0}, \mu_{k} R_{k})$,
and satisfy 
\begin{equation}
d(\theta_{i}^{k}, \theta_{j}^{k}) \geq\delta_{k} R_{k}
\end{equation}

where $\delta_{k} := \frac{1}{Ck}$, and $\mu_{k} := \frac{C}{k}$.

The points $\gamma_{k} \theta_{i}^{k}$ satisfy
\begin{equation}
d(\gamma_{k} \theta_{i}^{k}, \gamma_{k} \theta_{j}^{k}) \geq \delta,
\end{equation}

and $\frac{\epsilon_{k}}{\delta_{k}^{2}} = C^{2}$ is bounded.

Let us consider a point $\theta \in \mathcal{L}$. We want to approximate
it by a point of $\gamma_{k} D_{k}$. 

We can write $\theta = \gamma_{k} \theta_{k}$, for some 
$\theta_{k} \in \Lambda_{C'}$. There are two cases.

{\it Case 1). } For infinitely many indices $k$,
$\theta_{k} \in B_{\mathcal{L}}(\theta_{0}, R_{k})$.
We work in that case for these indices $k$, thus 
there are points $\theta'_{k} \in D_{k}\cap B_{\mathcal{L}}(\theta_{0}, R_{k})$,
with $d(\theta_{k},\theta'_{k}) \leq \epsilon_{k} R_{k}.$

Since the distance between the $\theta_{i}^{k}$'s is bounded below by
$\delta_{k} R_{k}$, we can find at least two of them which we call
$a_{k}$ and $b_{k}$, such that 

$$
d(\theta_{k}, b_{k}) \geq \frac{\delta_{k} R_{k}}{2}
$$

and,

\begin{equation}
d(\theta'_{k} , a_{k}) \geq \frac{\delta_{k} R_{k}}{2}.
\end{equation}

As $C'$ is contained in the cocompact group $\Gamma$, it
acts in a quasi-M\"{o}bius way on $(\partial \tilde{X},d)$

thus,

\begin{equation}
\frac{d(\gamma_{k} \theta'_{k}, \gamma_{k} \theta_{k}) d(\gamma_{k} a_{k}, \gamma_{k} b_{k})}
{d(\gamma_{k} \theta'_{k}, \gamma_{k} b_{k}) d(\gamma_{k} a_{k}, \gamma_{k} \theta_{k})}
\leq
\eta(\frac{d(\theta'_{k}, \theta_{k}) d(a_{k}, b_{k})}{ d(\theta'_{k}, b_{k}) d(\theta_{k}, a_{k})})
\end{equation}

for some homeomorphism $\eta : [0, \infty) \to [0, \infty)$.
This implies 

\begin{equation}
d(\gamma_{k} \theta'_{k}, \gamma_{k} \theta_{k})
\leq  \frac{(diam \mathcal{L}) ^{2}\eta(8\epsilon_{k} \mu_{k}/\delta_{k}^{2})}{\delta},
\end{equation}

therefore 
$d(\gamma_{k} \theta'_{k}, \gamma_{k} \theta_{k})$ tends to zero as 
$k$ tends to infinity.  

{\it Case 2) .} For all but finitely many indices $k$,
$\theta_{k} \notin B_{\mathcal{L}}(\theta_{0}, R_{k})$.

We work with these indices $k$ such that 
$\theta_{k} \notin B_{\mathcal{L}}(\theta_{0}, R_{k})$.

We know that $\epsilon_{k} / \delta_{k}^{2}$ is bounded above independantly of $k$,
and by assumption, $\delta_{k} \leq 2\mu_{k}$.

We claim that there exist $\xi_{k} \in  D_{k}$ 
and a positive constant $c_{0}$ such that for all $k$,

\begin{equation}
\frac{d(\xi_{k},\theta_{0})}{R_{k}} \geq c_{0}
\end{equation}
let us prove the claim.

On one hand, as $(\partial \tilde{X},d)$ is uniformly perfect, 
and so is it's weak tangent $(S,\delta)$ because 
the one point compactification $(\hat{S}, \hat{\delta})$
of  $(S,\delta)$
is quasi-M\"{o}bius homeomorphic to  $(\partial \tilde{X},d)$,
therefore there exist
a constant $C_{0} \in [0,1)$ such 
that for every $x \in  S$ and
$0 < R < diam  S$, we have 

\begin{equation}
\bar{B}_{(S,\delta)}(0,R) - B_{(S,\delta)}(0,C_{0}R)\neq \emptyset. 
\end{equation}

On the other hand, $(\mathcal{L}, \lambda_{k}d,\theta_{0})$ converges to
$(S, \delta,0)$.
After reindexing the sequence $\{ \lambda_{k} \}$, we have for each 
$\epsilon >0$ a map $g_{k} : B_{\lambda_{k} \mathcal{L}}(\theta_{0},k) \to S$ such that

\noindent
(i) $g_{k}(\theta_{0}) = 0$,

\noindent
for any two points $\theta$ and $\theta'$ in $B_{\lambda_{k} \mathcal{L}}(\theta_{0},k)$,

\noindent
(ii) $|\delta(g_{k}(\theta),g_{k}(\theta')) - \lambda_{k} d(\theta, \theta')| \leq \epsilon $,

\noindent
(iii) the $\epsilon$-neighborhood of $g_{k}(B_{\lambda_{k} \mathcal{L}}(\theta_{0},k))$
contains $B_{(S,\delta)}(0,k-\epsilon)$.

By (iii), we have

\begin{equation}
\bar{B}_{(S,\delta)}(0,k-\epsilon) \subset \mathcal{U}_{\epsilon}^{(S,\delta)}
g_{k}(\bar{B}_{\lambda_{k} \mathcal{L}}(\theta_{0},k)).
\end{equation}

By (6.13) there exist 
$y_{k} \in \bar{B}_{(S,\delta)}(0,k-\epsilon) - B_{(S,\delta)}(0,C_{0}(k-\epsilon))$, and by 
(6.14) there exist $\xi'_{k} \in \bar{B}_{\lambda_{k} \mathcal{L}}(\theta_{0},k)$ such that

\begin{equation}
\delta (y_{k}, g_{k}(\xi'_{k})) \leq \epsilon.
\end{equation}

We now evaluate $d(y_{k}, g_{k}(\xi'_{k}))$.
By the above properties (i), (ii), (6.15) and the triangle inequality we have

$$
\lambda_{k} d( \xi'_{k}, \theta_{0}) \geq \delta(g_{k}(\xi'_{k}), 0) -\epsilon
\geq \delta(y_{k},0) -\delta(y_{k},g_{k}(\xi'_{k})) -\epsilon
$$

\begin{equation}
\geq C_{0} (k-\epsilon) -2\epsilon.
\end{equation}

As $D_{k}$ is 
$\epsilon_{k} R_{k}$-dense in $B_{(\mathcal{L},d)}(\theta_{0},k/\lambda_{k})$,
there exist  $\xi_{k} \in D_{k}$ such that
 $d(\xi_{k},\xi'_{k}) \leq \epsilon_{k} R_{k} = \frac{k \epsilon_{k}}{\lambda_{k}}$.

Let us denote $c_{0}= C_{0}/2$. For $k$ large enough we have 
$\frac{C_{0}(k-\epsilon) -2\epsilon - k\epsilon_{k}}{\lambda_{k}} \geq \frac{c_{0} k}{\lambda_{k}}$,
therefore by (6.16) we get

\begin{equation}
d(\xi_{k},\theta_{0}) \geq d(\xi'_{k},\theta_{0}) - d(\xi'_{k},\xi_{k})\geq c_{0}R_{k},
\end{equation}

which proves the claim.
 
We can assume that for $k$ large enough,
$\mu_{k} < c_{0}/2 < 1/2$.

We choose $a_{k}= \theta_{1}^{k}$ and $b_{k}= \theta_{2}^{k}$, and we get 

$$
\frac{d(\gamma_{k} \xi_{k}, \gamma_{k} \theta_{k})
 d(\gamma_{k} a_{k}, \gamma_{k} b_{k})}
{d(\gamma_{k} \xi_{k}, \gamma_{k} b_{k})
 d(\gamma_{k} a_{k}, \gamma_{k} \theta_{k})}
\leq
\eta(\frac{d(\xi_{k}, \theta_{k})
 d(a_{k}, b_{k})}{ d(\xi_{k}, b_{k}) d(\theta_{k}, a_{k})})
$$

$$
\leq \eta( \frac{4\mu_{k}d(\theta_{k},\theta_{0}^{k})}{(d(\theta_{k},\theta_{0}^{k}) - \mu_{k}R_{k}))(c_{0}-\mu_{k})})
$$

\begin{equation}
\leq \eta(16\mu_{k}/c_{0}).
\end{equation}

We get

$$
d(\gamma_{k} \theta_{k}, \gamma_{k} \xi_{k}) \leq 
(diam\Lambda_{C'})^{2}\eta(16\mu_{k}/c_{0}) /\delta.
$$
$\Box$

\end{document}